\documentclass[12pt, colorinlistoftodos]{amsart}

 \usepackage{a4wide}
\usepackage{amssymb}
\usepackage{amsthm}
\usepackage{amsmath}
\usepackage{amscd}
\usepackage{verbatim}
\usepackage{bm}
\usepackage[mathscr]{eucal}
\usepackage[all]{xy}
\usepackage{todonotes}
\usepackage[style=alphabetic, backend=bibtex,
doi=false,isbn=false,url=false,backref=true,
eprint=true, sorting=nyt, maxnames=99]{biblatex}
\usepackage{hyperref}
\usepackage{dsfont, stmaryrd}
\usepackage{mathrsfs}
\usepackage{tikz, tikz-cd, caption, tabu}
\usepackage{tkz-euclide}
\numberwithin{equation}{section}

\theoremstyle{plain}
\newtheorem{theorem}{Theorem}[section]
\newtheorem{corollary}[theorem]{Corollary}
\newtheorem{lemma}[theorem]{Lemma}
\newtheorem{proposition}[theorem]{Proposition}

\newtheorem{conjecture}[theorem]{Conjecture}

\newtheorem{thmx}{Theorem}

\newtheorem{corx}[thmx]{Corollary}

\newtheorem{conjx}[thmx]{Conjecture}

\theoremstyle{definition}
\newtheorem{definition}[theorem]{Definition}
\newtheorem{remark}[theorem]{Remark}
\newtheorem{example}[theorem]{Example}

\theoremstyle{remark}

\newcommand{\pmat}[4]{\begin{pmatrix}
                 #1 & #2\\
                 #3 & #4
\end{pmatrix}}
\newcommand{\smat}[4]{\left(\begin{smallmatrix}
                 #1 & #2\\
                 #3 & #4
\end{smallmatrix}\right)}

\newcommand{\lp}{\left (}
\newcommand{\rp}{\right )}

\newcommand{\Fc}{{\mathcal{F}}}

\newcommand{\Oc}{{\mathcal{O}}}
\newcommand{\Sc}{{\mathcal{S}}}
\newcommand{\Zb}{\mathbb{Z}}

\newcommand{\Qb}{\mathbb{Q}}
\newcommand{\SO}{{\mathrm{SO}}}
\newcommand{\SL}{{\mathrm{SL}}}
\newcommand{\GL}{{\mathrm{GL}}}
\newcommand{\sgn}{{\mathrm{sgn}}}

\newcommand{\df}{\mathfrak{d}}

\newcommand{\pf}{\mathfrak{p}}
\newcommand{\mf}{\mathfrak{m}}

\newcommand{\Nm}{{\mathrm{Nm}}}
\newcommand{\Ab}{\mathbb{A}}
\newcommand{\Db}{\mathbb{D}}

\newcommand{\Nb}{\mathbb{N}}
\newcommand{\Hb}{\mathbb{H}}
\newcommand{\Rb}{\mathbb{R}}
\newcommand{\Cb}{\mathbb{C}}

\newcommand{\Gm}{\mathbb{G}_m}
\newcommand{\ebf}{{\mathbf{e}}}

\newcommand{\tH}{{\tilde{H}}}

\renewcommand{\th}{\tilde{h}}

\newcommand{\tr}{\operatorname{Tr}}

\newcommand{\Gal}{\operatorname{Gal}}

\newcommand{\dd}{\mathrm{d}}

\newcommand{\tphi}{\tilde\phi}

\newcommand{\slf}{{\mathfrak{sl}}}

\newcommand{\vol}{\mathrm{vol}}

\newcommand{\tg}{\tilde{g}}

\newcommand{\tG}{\tilde{G}}

\newcommand{\cL}{\mathcal{L}}

\newcommand{\cha}{\mathrm{Char}}

\newcommand{\Cl}{\mathrm{Cl}}

\newcommand{\Tc}{\mathcal{T}}
\newcommand{\Cc}{\mathcal{C}}
\newcommand{\Nc}{\mathcal{N}}
\renewcommand{\aa}{A}
\newcommand{\fs}{f^\#}
\newcommand{\hf}{\mathfrak{h}}
\newcommand{\gf}{\mathfrak{g}}
\newcommand{\gl}{\mathfrak{gl}}

\newcommand{\supp}{\mathrm{supp}}
\AtBeginDocument{%
   \def\MR#1{}
}

\newcommand{\Q}{\mathbb{Q}}
\newcommand{\C}{\mathbb{C}}
\newcommand{\R}{\mathbb{R}}
\newcommand{\Z}{\mathbb{Z}}

\renewcommand{\P}{\mathbb{P}}

\newcommand{\Spec}{\operatorname{Spec}}
\newcommand{\Hom}{\operatorname{Hom}}
\newcommand{\Reg}{\operatorname{Reg}}
\newcommand{\Ad}{\operatorname{Ad}}
\newcommand{\Ind}{\operatorname{Ind}}

\newcommand{\iso}{\cong}
\renewcommand{\to}{\rightarrow}
\newcommand{\tor}{\mathrm{tor}}

\renewcommand{\O}{\mathcal{O}}

\begin{document}
\title{Motivic action conjecture for Doi--Naganuma lifts}

\author[A.~ Horawa]{Aleksander Horawa}
\address{Mathematical Institute\\University of Oxford\\Woodstock Road\\Oxford, OX2
	6GG, UK}
\email{horawa@maths.ox.ac.uk}

\author[Y.~Li]{Yingkun Li}
 \address{
Max Planck Institute for Mathematics,
    Vivatsgasse 7, 
    D--53111     Bonn,
    Germany}
\email{yingkun@mpim-bonn.mpg.de}

\subjclass[2020]{}
\thanks{
}
\begin{abstract}
  We prove the motivic action conjecture for the base change to real quadratic fields of weight one newforms with odd, squarefree level and solvable projective image.
\end{abstract}
\date{\today}
\maketitle

 
\section{Introduction}

The motivic action conjecture, introduced by Akshay Venkatesh and his collaborators~\cite{Prasanna_Venkatesh, Venkatesh:Derived_Hecke, Harris_Venkatesh}, predicts a precise relationship between the rationality of contributions of automorphic forms to cohomology and certain motivic cohomology groups. A concrete instance of this phenomenon~\cite{Horawa} is that logarithms of Stark units should be factors of rationality for coherent cohomology classes associated with weight one Hilbert modular forms. We prove this conjecture for certain Hilbert modular forms over real quadratic fields which arise as base change of weight one newforms.
Historically, such Hilbert modular forms in higher weight were studied by Doi and Naganuma using the converse theorem \cite{DN69}. Later, they were realized as theta lifts from $\SL_2$ to $\mathrm{O}(V)$ for a suitable quadratic space $V$ of signature $(2, 2)$ by Kudla \cite{Kudla78} (see Section \ref{subsec:theta-lifts}).

Let $F$ be a real quadratic field of
discriminant $D > 0$, and $f$ be a Hilbert modular form of parallel weight one for $F$. There is a line bundle $\omega$, the {\em Hodge bundle}, on the Hilbert modular surface $X$ defined over $\Q$, such that $f \in H^0(X, \omega) \otimes_\Q \Q(f)$, where $\Q(f)$ is the number field generated by the Hecke eigenvalues of $f$. Considering the higher cohomology of $\omega$, there are natural classes~\eqref{eqn:omega_f^j}:
$$\omega_f^1, \omega_f^2 \in H^1(X, \omega) \otimes_\Q \C$$
associated with $f$. We are interested in the rationality properties of $\omega_f^1, \omega_f^2$.

The Galois representation associated with $f$ is an odd Artin representation 
$$\varrho_f \colon \Gal(L/ F) \to \GL_2(\Q(f))$$ 
for some number field $L$. The associated conjugation action of $\Gal(L/F)$ on traceless $2 \times 2$ matrices over $\Q(f)$ gives rise to the representation 
$$\Ad^0 \varrho_f \colon \Gal(L/F) \to \GL_3(\Q(f)).$$
Stark conjectured~\cite{Stark-II} that there are units $u_{11}, u_{12}, u_{21}, u_{22} \in \O_L^\times \otimes_\Z \Q(f)$ such that the second derivative of the Artin $L$-function at $s = 0$ is expressed in terms of their logarithms:
$$L^{(2)}(\Ad^0 \varrho_f, 0) \sim_{\Q(f)^\times}  \det R_f \qquad \text{for } R_f := \begin{pmatrix}
	\log|u_{11}| & \log|u_{12}| \\
	\log|u_{21}| & \log |u_{22}|
\end{pmatrix},$$
where $a \sim_{\Q(f)^\times} b$ if there exists $\lambda \in \Q(f)^\times$ such that $a = \lambda b$. He proved it when $\Ad^0 \varrho_f$ has rational traces (loc.\ cit.). 

The motivic\footnote{These Stark unit groups are simple examples of {\em motivic cohomology groups} which are predicted to play a similar role in a more general story~\cite{Prasanna_Venkatesh, Horawa_Prasanna}.} action conjecture for Hilbert modular forms relates the rationality of $\omega_f^1, \omega_f^2$ to the regulator matrix $R_f$.

\begin{conjx}[{\cite{Horawa}}]\label{conjA}
	The following coherent cohomology classes associated with $f$ are rational:
$$
	\frac{\log|u_{22}| \omega_f^1 - \log|u_{12}| \omega_f^2}{\det R_f}, \frac{-\log|u_{21}| \omega_f^1  + \log|u_{11}| \omega_f^2}{\det R_f} \in H^1(X, \omega) \otimes_\Q \Q(f).
$$
\end{conjx}

\begin{thmx}[Theorem~\ref{thm:motivic_action}]\label{thmB}
  Let $f$ be the quadratic base change of a weight one newform $f_0 \in S_1(N, \chi_0)$.
  Suppose $D, N$ are odd, squarefree, co-prime to each other, and $f_0$ has projective image not isomorphic to $A_5$. Then Conjecture~\ref{conjA} is true. 
	Explicitly, there are units $u_{f_0}$ associated with $\Ad^0 \varrho_{f_0}$ and $u_{f_0}^F$ associated with the quadratic twist $\Ad^0 \varrho_{f_0} \otimes \chi_F$ such that:
	$$  \frac{\omega_f^1 + \omega_f^2}{\log|u_{f_0}|}, \frac{\omega_f^1 - \omega_f^2}{\log|u_{f_0}^F|} \in H^1(X, \omega) \otimes \Q(f).$$
\end{thmx}

\begin{remark}
	\leavevmode 
  \begin{enumerate}
  	\item The conditions that $D, N$ are odd and squarefree are technical, and can be removed with some additional calculations of the Fourier expansion in Section \ref{subsec:schwartz}. The condition that $D, N$ are co-prime is a bit more subtle to remove.
  	\item   If $f_0$ has projective image $A_5$, then Theorem \ref{thmB} also holds assuming Stark's conjecture for $\mathrm{Ad}^0\varrho_{f_0}$ and $\Ad^0 \varrho_{f_0} \otimes \chi_{F}$.
  \end{enumerate}
\end{remark}

For Hilbert modular forms $f$ of higher weight, Michael Harris~\cite{Harris_Periodinvariants} defined period invariants $\nu_f^1, \nu_f^2 \in \C$ by requiring that $\omega_f^j/\nu_f^j$ be rational in the appropriate coherent cohomology group. These give periods of rationality for Rankin--Selberg, triple product, and Asai L-functions of Hilbert modular forms; for example, the coherent cohomology classes $\omega_f^j/\nu_f^j$ were recently interpolated $p$-adically to construct an Asai $p$-adic $L$-function for Hilbert modular forms~\cite{Grossi_Loeffler_Zerbes}.

This raised the natural question: are any multiples of $\omega_f^1, \omega_f^2$ rational in coherent cohomology for a weight one Hilbert modular form $f$? Theorem~\ref{thmB} finally settles this question.

\begin{corx}[Corollary~\ref{cor:Harris_periods}]
	Under the assumptions of Theorem~\ref{thmB}, no multiple of $\omega_f^1$ or $\omega_f^2$ belongs to $H^1(X ,\omega) \otimes_{\Q} \overline \Q$.
\end{corx}

We refer to Section~\ref{subsec:period_invariants} for more details.

\subsection{Sketch of proof}

The classes $\omega_f^1, \omega_f^2$ give a complex basis of the $f$-isotypic component $H^1(X, \omega)_f \otimes \C$ and we are interested in the rational structure $H^1(X, \omega)_f$ on this vector space. Using the non-trivial element in $\Gal(F/\Q)$, we can define an involution $s^\ast$ on $H^1(X, \omega)$ which preserves the $f$-isotypic components if $f$ comes from base change from $\Q$. This breaks $H^1(X, \omega)_f$ into two 1-dimensional spaces $H^1(X, \omega)_f^\pm$ where the superscript $\pm$ means $s^\ast = \pm 1$. Over $\C$, $s^\ast \omega_f^1 = \omega_f^2$, and hence our task is to show that:
\begin{align}
	\eta^+ & := \frac{\omega_f^1 + \omega_f^2}{\log |u_{f_0}|} \in H^1(X, \omega)_f^+ \\
	\eta^- & := \frac{\omega_f^1 - \omega_f^2}{\log |u_{f_0}^F|} \in H^1(X, \omega)_f^-.
\end{align}
Serre duality pairs $H^1(X, \omega)_f^+$ with $H^1(X, \omega)_f^-$ and we check that $\eta^+$ pairs rationally with~$\eta^-$. This shows that $\eta^+$ is rational if and only if $\eta^-$ is, and hence reduces Theorem~\ref{thmB} to proving that $\eta^-$ is rational.

We construct the quadratic base change $f$ of a weight one modular form $f_0$ explicitly as a theta lift. We then use this explicit construction to give a rationality criterion for $H^1(X, \omega)_f^-$ by integrating the cohomology class against a linear combination of Hirzebruch--Zagier divisors.
The statement below implies that $\eta^- \in H^1(X, \omega)_f^-$.

\begin{thmx}[see Theorem~\ref{thm:unfold2} and Corollary~\ref{cor:Cc}]\label{thmD}
  There exists a linear combination of cycle integrals along Hirzebruch-Zagier divisors that define an associated functional
  $$
  \mathcal C \colon H^1(X, \omega)_f \to \Q(f),
  $$ 
  such that $\mathcal C_\C \colon H^1(X, \omega)_f \otimes \C \to \C$ satisfies $\mathcal C_\C(\eta^-) \in \Q(f)^\times$ and $\mathcal C_C(\eta^+) = 0$.
\end{thmx}

The linear combination above is naturally phrased in terms of weighted cycle as in \cite{Kudla97}. 
The associated functional appears in the Fourier expansion of the theta lift of $f$ to a scalar multiple $r_0$ of $f_0$ through a standard  unfolding process.
The key point  is to compute this scalar $r_0$, which depends on the Schwartz function $\varphi$ defining the theta function.
In fact, there is a family of functionals, whose generating series maps $f$ to $r_0 f_0$. 

It is not difficult to show that $r_0$ is a $\overline\Qb^\times$-multiple of the ratio $\frac{\langle f, f \rangle}{\langle f_0, f_0 \rangle}$. However, it is rather tricky to show that the associated functional $\Cc$ takes values in $\Qb(f)$. 
The functional $\Cc$ appearing in the Fourier expansion involves certain linear combination of connected components of the weighted cycle, which are then defined over certain extensions of $\Qb$.
It turns out that the Galois action on the components and the coefficients compensate each other, giving us $\Cc(\eta^-) \in \Qb(f)$.  
This phenomenon seems to only appear when $f_0$ has non-trivial nebentypus, in particular for odd weight. 

\begin{remark}
	Theorem~\ref{thmD} is the analogue of the following period relation proved by Oda~\cite[Section 16]{Oda_book}. Suppose that $F$ has a fundamental unit $\epsilon \in \O_F^\times$ with $\epsilon_1 > 0$, $\epsilon_2 < 0$ (i.e.\ $N_{F/\Q}(\epsilon) = -1$), and $f$ be a Hilbert modular form of parallel weight $k \geq 2$ and level one. Then Oda proved:
	$$\int_{\SL_2(\Z) \backslash \mathbb H} f(\epsilon_1 z, \epsilon_2 \overline z) y^k \frac{dx dy}{y^2} = \begin{cases}
		0 & \text{$f$ is not a Doi--Naganuma lift} \\
		c \frac{\langle f, f \rangle}{\langle f_0, f_0 \rangle} & \text{$f$ is the Doi--Naganuma lift of $f_0$} \\
	\end{cases}$$
	for some $c \in \Q^\times$. Note that the Hirzebruch--Zagier divisor here is simply the diagonally embedded modular curve, and 
	$$\omega_f^2 = \left[ f(\epsilon_1 z_1, \epsilon_2 \overline z_2) y^2 \frac{dx_2 dy_2}{y_2^2} \right].$$
\end{remark}

\begin{remark}
	The first-named author had previously obtained numerical evidence for the conjecture by restricting to diagonally embedded modular curves~\cite[Section 6]{Horawa}. In that case, the restriction was non-zero on $\eta^+$ and zero on $\eta^-$, and one could numerically verify the first restriction is rational (with many digits of accuracy). Somewhat surprisingly, the Hirzebruch--Zagier divisors studied in this paper have the opposite property --- the restrictions are non-zero on $\eta^-$ and zero of $\eta^+$. As a roundabout consequence, we do obtain the rationality of $\eta^+$, but our proof is genuinely different than the numerical verification.
\end{remark}

\subsection{Potential generalizations}

We end the introduction by explaining what we expect to happen for general quadratic base change. Suppose $F/F_0$ is a quadratic extension of totally real fields. Starting with a Hilbert modular form $f_0$ of parallel weight one for $F_0$, one can again construct its base change $f$ to $F$ as a theta lift. For each subset of the places of $F$, $J \subseteq \{1, \ldots, d\}$, we have a cohomology class:
$$\omega_f^J \in H^{|J|}(X, \omega)_f \otimes \C,$$
where $X$ is a Hilbert modular variety for $F$ and $\omega$ is the Hodge bundle on $X$. The unit group $U_f$ has rank $d = 2d_0$, and we expect $j$th exterior power of the $d \times d$ Stark regulator matrix $R_f$ to predict the rationality in $H^j(X_\C, \omega)_f$ \cite[Conjecture 4.17]{Horawa}.

The methods of this paper should give results about rationality in middle degree coherent cohomology, i.e.\ $H^{d_0}(X_\C, \omega)_f$. The non-trivial element $\sigma$ of $\Gal(F/F_0)$ acts on the set of infinite places of $F$, and it is natural to consider the subspace
$$V_0^\C = \mathrm{span}\{ \omega_f^J \ | \ |J| = d_0, J \cap J^\sigma = \emptyset  \} \subseteq H^{d_0}(X, \omega)_f \otimes \C$$
of dimension $2^{d_0}$ within the $\binom{2d_0}{d_0}$-dimensional space. Restriction of classes in $H^{d_0}(X_\C, \omega)_f$ to Shimura subvarieties associated with $F_0$ can only be non-zero for elements in $V_0^\C$. The motivic action conjecture predicts that this subspace is rational, i.e.\ $V_0^\C = (V_0)_\C$ for a rational subspace $V_0 \subseteq H^{d_0}(X, \omega)_f$. The swap map $s^\ast$ preserves $V_0$, breaking it up into $2^{d_0-1}$-dimensional isotypic components $V_0 = V_0^+ \oplus V_0^-$, paired with each other under Serre duality. Finally, the analogue of Theorem~\ref{thmD} would give rise to $2^{d_0-1}$ functionals on $V_0^-$ which take rational values on the expected basis of $V_0$. More concretely, assuming the subspace $V_0 \subseteq H^{d_0}(X, \omega)_f$ is rational, we would hope to prove that
$$ \frac{\omega_f^J + \omega_f^{J^\sigma}}{\det R_{f_0}},   \frac{\omega_f^J - \omega_f^{J^\sigma}}{\det R_{f_0}^F} \in H^{d_0}(X, \omega)_f$$
for any $J$ with $|J| = d_0$ and $J \cap J^\sigma = \emptyset$, where $R_{f_0}$ and $R_{f_0}^F$ are the $d_0 \times d_0$ regulator matrices for $\Ad^0 \varrho_0$ and $\Ad^0 \varrho_0 \otimes \chi_{F/F_0}$, respectively. We hope to pursue this generalization in future work.

The core idea behind our approach is that if an automorphic form is obtained by theta lifting from a smaller group, then one can give a rationality criterion for coherent cohomology by considering the theta lift back. The same strategy also led to the proof of the motivic action conjecture for Siegel modular forms in special cases~\cite{Horawa_Prasanna}. It would be interesting to explore other instances of this phenomenon.

\subsection*{Acknowledgements}

We would like to thank Gerard van der Geer, James Newton, Kartik Prasanna, and Akshay Venkatesh for their interest and helpful discussions.

AH was supported by NSF grant DMS-2001293 and UK Research and Innovation grant MR/V021931/1.
YL was supported by the Deutsche Forschungsgemeinschaft (DFG) through the Collaborative Research Centre TRR 326 ``Geometry and Arithmetic of Uniformized Structures'' (project number 444845124) and the Heisenberg Program ``Arithmetic of real-analytic automorphic forms''(project number 539345613).
This work was completed while YL was visiting MPIM Bonn. He thanks the institute for its friendly working environment. 

For the purpose of Open Access, the authors have applied a CC BY public copyright licence to any Author Accepted Manuscript (AAM) version arising from this submission.

\setcounter{tocdepth}{1}
\tableofcontents

\section{Preliminaries}\label{sec:prelim}
In this section, we recall preliminaries on automorphic forms, and their theta lifts. 

\subsection{Notation}
\label{subsec:notations}
For a ring $R$ and any $N \in \Nb$, we denote $R_N :=  R \otimes_\Zb\prod_{p \mid N \text{ prime}}\Zb_p \subset \hat R := R \otimes \hat \Zb$. 
Similarly for $a \in R$, we write $a_N := a \otimes 1 \in R_N \subset \hat R$. 
Throughout the paper,~$F$ will be a totally real field of degree $d$ with different $\df$, ring of integers $\Oc$, and adeles $\Ab_F$. Apart from Section \ref{subsec:Hmf}, $d$ will be 2, i.e.\ $F$ will be real quadratic. In that case, let $D > 0$ for the discriminant of $F$, which we assume is odd throughout for simplicity.

For a compact open subgroup $U \subset \hat\Oc^\times$, we define the associated narrow class group 
\begin{equation}
  \label{eq:Cl+}
  \Cl^+_F(U) :=
   \Ab_F^\times /F^\times U \Rb_{>0}^d.
 \end{equation}
 We omit $F$ and $+$ in then notation when $F = \Qb$. For an algebraic group $G$ defined over $\Qb$, we use $[G]$ to denote the quotient $G(\Qb)\backslash G(\Ab)$. 
 
For a quadratic extension $K/\Q$, we will write $\chi_K$ for the associated quadratic character of $\Q$.
 
\subsection{Adelic Hilbert modular forms}
\label{subsec:Hmf}
We first recall adelic Hilbert modular form (see \cite[\textsection 1, 2]{Shimura78-dmj}, \cite{Garrett-book}). 
Denote $G_F := \mathrm{Res}_{F/\Qb}\GL_2$ with center $Z_F = \mathrm{Res}_{F/\Qb}\Gm \subset G_F$. 
We omit $F$ from the notation when it is $\Qb$. 
Denote $K_\infty = \SO_2(\Rb)^d \subset G_F(\Rb) $ the maximal compact. 
For open compact $K \subset G_F(\hat\Qb)$, we have the decomposition
\begin{equation}
  \label{eq:double-coset}
  G_F(\Ab) = \bigsqcup_{\xi \in \Xi} G_F(\Qb) \xi K G_F(\Rb)^+
\end{equation}
for a finite subset $\Xi \subset G_F(\hat\Qb)$ satisfying
\begin{equation}
  \label{eq:AFx-decomp}
 \Ab_F^\times = \bigsqcup_{\xi \in \Xi} F^\times \det(\xi) \det(K) \Rb_{>0}^d,
\end{equation}
where $G_F(\Rb)^+ \cong (\GL_2(\Rb)^+)^d$ is the subgroup of $G_F(\Rb) \cong \GL_2(\Rb)^d$ with positive determinant.
This decomposes $G_F(\Qb)\backslash G_F(\Ab)/K$ into cosets
\begin{equation}
  \label{eq:double-coset1}
    [G_F]/K := G_F(\Qb)\backslash G_F(\Ab)/K = \bigsqcup_\xi \Gamma_\xi \backslash G_F(\Rb)^+
\end{equation}
with $\Gamma_\xi := G_F(\Qb)^+ \cap \xi K \xi^{-1} $, $G_F(\Qb)^+ := G_F(\Qb) \cap G_F(\Rb)^+$.
The class group $\Cl^+_F(\det(K))$ acts simply transitively on $\Xi$.
 We will choose $\Xi$ to contain the identity, and use it as a base point to identify the sets $\Xi \cong \Cl^+_F(\det(K))$. 
 If $\det(K) = \hat\Oc^\times$, then $\Cl^+_F(\det(K))$ is the narrow class group $\Cl^+_F$ of $F$. 
Similarly, if we denote  $G_{1, F} := \mathrm{Res}_{F/\Qb}\SL_2 \subset G_F$, then there is a finite set $\Xi' \subset \Xi$ such that
\begin{equation}
  \label{eq:double-coset2}
  G_F(\Ab) = \bigsqcup_{\xi' \in \Xi'} Z_F(\Ab) G_F(\Qb)  G_{1, F}(\Ab) \xi' K
\end{equation}
and $\Ab_F^\times = \sqcup_{\xi \in \Xi'} (\Ab_F^\times)^2 F^\times \det(\xi') \det(K)$.
The subset of $\Cl^+_F(\det(K))$ corresponding to $\Xi' \subset \Xi$ consists of representatives of
$\Cl^+_F(\det(K))/\Cl^+_F(\det(K))^2$.

On $\Ab_F$, let $db$ be the Haar measure self-dual with respect to the additive character $\psi_F := \psi \circ \Nm_{F/\Qb}$, where
 $\psi = \otimes_{p \le \infty} \psi_p$ is the additive character of $\Qb\backslash \Ab$ normalized with $\psi_\infty(x) = \ebf(x) := e^{2\pi i x}$. 
 It is the product of local measures $db_v$ over finite places $v$ of $F$. With respect to $db_v$, the valuation ring $\Oc_v$ of $F_v$, the completion $F$ at $v$, has volume
 $ |\df|_v^{1/2}$.
 On $\Ab^\times_F$, we normalize the Haar measure $d^\times a = \prod_{v \le \infty} d^\times a_v$ such that $d^\times a_v =  m_v\frac{da_v}{|a|_v}$ with $m_v = (1-1/q_v)^{-1}$ and $q_v = |\Oc_v/\varpi_v|$ for $v < \infty$, and $d^\times a_\infty = \frac{da_\infty}{|a|_\infty}$.
Then $\Oc_v^\times$ has volume 1 with respect to $d^\times a_v$.

Let $dg = \prod_{v \le \infty} dg_v$ be the Haar measure on $G_F(\Ab) = G(\Ab_F) = \otimes_{v \le \infty} G(F_v)$ normalized such that the subgroup $\SL_2(\Oc_v) \subset G(F_v)$ has volume 1 with respect to $dg_v$ for all finite places $v$ of $F$.
At the infinite places, we have $dg_\infty = \prod_{1 \le j \le d} d\mu(z_j)\frac{d \theta_j}{2\pi}$ in the coordinate $g = (g_{z_j} k_{\theta_j}) \in \SL_2(\Rb)^d \cong G_{1, F}(\Rb) $, where for $z = x + i y \in \Hb, g_z =
n(x)m(\sqrt{y}) = \smat{1}{x}{}{1} \smat{\sqrt{y}}{}{}{1/\sqrt{y}}, 
d\mu(z) = \frac{dx dy}{y^2},$ and $k_\theta = \smat{\cos \theta}{\sin \theta}{-\sin \theta}{\cos\theta} \in \SO_2(\Rb)$. 

Let $k_1, \dots, k_d$ be integers satisfying $k_j \equiv k_1 \bmod 2$.
An automorphic form of weight $(k_1, \dots, k_d)$, level $K$ and central character $\chi: F^\times \backslash \Ab_F \to \Cb$ is a function $f: [G_F]/K \to \Cb$ right $K_\infty$-equivariant with respect to the character $(\theta_1, \dots, \theta_d) \mapsto e^{i \sum_j k_j \theta_j}$ and satisfying
\begin{equation}
  \label{eq:center}
  \rho(z) f = \chi(z) f
\end{equation}
for all $z \in Z_F(\Ab)$, 
where we have a left action on $f$ given by
\begin{equation}
  \label{eq:rho-g0}
  (\rho(g) f)(g') := f(g g'). 
\end{equation}
for any $g' \in G_F(\Ab)$.
This right action induces an action of the Lie algebra $\gl_2^d$ of $G_F(\Rb) \cong \GL_2(\Rb)^d$.
For $1 \le j \le d$, denote $L_j, R_j \in \slf_{2, \Cb}^d \subset \gl_{2, \Cb}^d$ the element whose $j$-th component is given by
\begin{equation}
  \label{eq:LR}
  L =  \frac12 \pmat{1}{-i}{-i}{-1},~ 
R =   \frac12 \pmat{1}{i}{i}{-1}
\end{equation}
respectively, and other components are 0. 
Their actions on $f$ correspond to the lowering and raising operator in the $j$-th variable.
We say $f$ is holomorphic, resp.\ anti-holomorphic, in the $j$-th variable if $L_j f$, resp.\ $R_j f$, vanishes.
We call it (anti-)holomorphic if it is (anti-)holomorphic for all $1 \le j \le d$.

If $\chi$ is unitary, then $|f|$ is automorphic of weight 0, level $K$ and trivial central character.
Furthermore, if $|f|$ is in $L^2([G_F]/Z_F(\Ab))$, e.g.\ $f$ is a cusp form, then we define the Petersson norm of $f$ to be
\begin{equation}
  \label{eq:Pet}
  \|f\|_{\mathrm{Pet}}^2 :=
  \int_{G_F(\Qb)\backslash G_F(\Ab)/Z_F(\Ab)} |f(g)|^2 dg,~ 
  \|f\|_{1, \mathrm{Pet}}^2 :=
  \int_{G_{1, F}(\Qb)\backslash G_{1,F}(\Ab)} |f(g)|^2 dg,   
\end{equation}
Using the decomposition \eqref{eq:double-coset2}, we can rewrite
\begin{equation}
  \label{eq:Pet-sum}
  \|f\|_{\mathrm{Pet}}^2 :=
\sum_{\xi' \in \Xi'}  \|\rho(\xi')f\|_{1, \mathrm{Pet}}^2.
\end{equation}
In particular if $\Xi'$ has size 1, i.e.\ $\Cl^+_F(\det(K))$ has odd cardinality, then $  \|f\|_{\mathrm{Pet}}^2 :=
\|f\|_{1, \mathrm{Pet}}^2$.

The decomposition \eqref{eq:double-coset1} shows that $f$ corresponds to a tuple of classical Hilbert modular forms on $\prod_{\xi \in \Xi} \Gamma_\xi \backslash \Hb^d$.
We denote $f^\flat(z) = (f^\flat_\xi(z))_{\xi \in \Xi}$ the classical Hilbert modular form associated to $f$ for $z = (z_j)_{1 \le j \le d} \in \Hb^d$.
If $f$ is a holomorphic cusp form, then it has the Fourier expansion
\begin{equation}
  \label{eq:fb-FE}
  f^\flat_\xi(z)
  =
  \sum_{\nu \in F,~ \nu \gg 0}
c_{\nu, \xi}(  f)  \ebf(\tr(\nu z)). 
\end{equation}
The function $  f^\flat_\xi$ satisfies the defining property
\begin{equation}
  \label{eq:adele-classical}
  f(\gamma g) = f_\xi^\flat(g_\infty \cdot i) \prod_{1 \le j \le d} \mathrm{j}(g_{\infty, j}, i)^{-k_j},~ \mathrm{j}(\smat{a}{b}{c}d, z) = (cz + d) (ad-bc)^{-1/2}
\end{equation}
for $\gamma \in G_F(\Qb), g = (g_f, g_\infty) \in \xi K G_F(\Rb)^+$.
We omit $\xi$ from the notation when it corresponds to the different ideal $\df$ and $\Xi$ has size 1. 

Similarly for a classical modular form $f_0$, we use $f_0^\#$ to denote its adelization. 
For example when $F = \Qb$ and $f_0 \in S_{k}(N, \chi_0)$, then $\fs_0$
satisfies (see \cite[Prop.\ 1.4]{Kudla-book})
\begin{equation}
  \label{eq:f0-adel}
  f_0^\#(\gamma g \kappa) = \chi_0(\kappa)f_0^\#(g)
\end{equation}
for all $\gamma \in \GL_2(\Qb), g \in \GL_2(\Ab)$ and
\begin{equation}
  \label{eq:K0N}
  \kappa \in K_0(N)
  :=\{\gamma \in \GL_2( \hat\Zb): \gamma \equiv \smat{*}{*}{0}{*} \bmod N\}.
\end{equation}
Here we also used $\chi_0 = \otimes_{p \le \infty} \chi_{0, p}: \Qb^\times\backslash\Ab^\times \to \Cb^\times$  to denote the associated idelic character, and define $\chi_0(\kappa) := \chi_0(a)$ for $\kappa = \smat{a}{*}{*}{*}$.
It is easy to check that
\begin{equation}
  \label{eq:cmf0}
\sqrt{v}^{k}  c_m(f_0) e^{-2\pi m v} = \int_{\Qb\backslash \Ab} \fs_0(n(b) m(\sqrt{v})) \psi(-mb) db
\end{equation}
for any $m \in \Qb$. 



\subsection{Connected components of adelic Hilbert modular surfaces}

Let $F/\Qb$ be a real quadratic field.
In this section, we recall the adelic Hilbert modular surfaces associated with a suitable subgroup $ H \subset G_F$.
They are call geometric in the literature as the associated Shimura variety has good moduli interpretation. 

Denote $ \tG = \GL_2$ and $\tH = G_F \times \Gm$. 
Define the embedding
\begin{equation}
  \label{eq:emb-d}
  \dd: \mathbb{G}_m \to \GL_2,~ \alpha \mapsto \pmat{1}{}{}{\alpha}. 
\end{equation}
On $\tG \times \tH$, we have the character
\begin{equation}
  \label{eq:similitude}
\nu: \tG \times \tH \to \mathbb{G}_m,~ 
  \nu(\tg, (h, a)) = \det(g)\Nm(\det(h)/a),
\end{equation}
and are interested in the subgroup
\begin{equation}
  \label{eq:simi}
  R :=  \{(\tg, \tilde h) \in
\tG \times \tH
  : \nu(\tg) = \nu(\tilde h)\}.
\end{equation}
It contains $G \times H$, where $ G := \ker \nu \cap G = \SL_2$
and \begin{equation}
  \label{eq:tH1}
  \begin{split}
    H
    &:=
G_F \times_{\mathrm{Res}_{F/\Qb} \Gm} \Gm
      \subset  \ker \nu \cap \tH.
  \end{split}
\end{equation}
Through projection to the first factor, we view $H$ as a subgroup of $G_F$. 
Note $R\cong G \rtimes \tH$ via the map $(\tg, h) \mapsto (g, h)$, where
\begin{equation}
  \label{eq:g1}
 g := \tg \cdot \dd(\nu(\tg)^{-1}). 
\end{equation}
In addition, we also denote 
\begin{equation}
  \label{eq:H1}
 H_1 := \mathrm{Res}_{F/\Qb}\SL_2 \subset H \subset G_F .
\end{equation}

For an open compact subgroup $K \subset H(\hat\Qb)$, we can form the Shimura variety $X_K$. 
It is defined over $\Qb$ and its $\Cb$-points are
\begin{equation}
  \label{eq:XK}
  X_K(\Cb) =
  H(\Qb) \backslash (\Db \times  H(\hat\Qb)/K),
\end{equation}
where $\Db = \Hb^2 \cup (\Hb^-)^2$ is the symmetric space associated to $H/Z_H$ and $Z_H \subset H$ the center.
As in \eqref{eq:double-coset2}, we can write
\begin{equation}
  \label{eq:double-coset3}
  H(\Ab) = \bigsqcup_{\delta \in \Delta} H(\Qb) H_1(\Ab)\delta K Z_H(\Rb)
  \end{equation}
  for a finite subset $\Delta \subset H(\hat\Qb)$ satisfying
  $$
\hat\Qb^\times =   \bigsqcup_{\delta \in \Delta}
  \Qb^\times \det(\delta) \det(K).
  $$
We identify the sets $\Delta \cong \Cl(\det(K)) = \hat\Qb^\times/\Qb^\times\det(K)$, and let the latter acts the former by left multiplication. 
  This leads to a decomposition into components
  \begin{equation}
    \label{eq:X1K}
    X_{K} = \bigsqcup_{\delta \in \Delta}
    X_{1, K,{\delta}},
  \end{equation}
  where
  $X_{1, K, \delta}$ is the geometrically irreducible Shimura variety associated to $H_1$ and open compact subgroup
\footnote{In particular, $K_1 = K \cap H_1(\hat\Qb)$.}
  $K_{\delta} := \delta K \delta^{-1}  \cap H_1(\hat\Qb) \subset H_1(\hat\Qb)$.
  It is important to know that this decomposition is defined over the abelian extension $E_K/\Qb$ determined by the reciprocity map
  \begin{equation}
    \label{eq:E}
\Cl(\det(K))
    \cong \Gal(E_K/\Qb),~ a \mapsto \sigma_a^{-1}.
  \end{equation}
The Galois group $\Gal(E_K/\Qb)$ acts on the set $\{ X_{1, K, {\delta}}: \delta \in \Delta\}$ through \cite[(1.10), (1.11)]{Kudla97}
  \begin{equation}
    \label{eq:E-act}
    X_{1, K,{\delta}}^{\sigma}
    \cong     X_{1, K,{\sigma^{}\cdot \delta }}.
  \end{equation}
  So decomposition \eqref{eq:X1K} can be rewritten as
  \begin{equation}
    \label{eq:X1K-Gal}
    X_{K} \cong \bigsqcup_{\sigma \in \Gal(E_K/\Qb)}
    X^\sigma_{1, K_1},~
    X_{1, K_1} := X_{1, K, 1}.
  \end{equation}

  \subsection{Weil representation}

Let $F$ be a real quadratic field of discriminant $D > 0$ as in the previous section.
We recall the Weil representation $\omega = \omega_\psi$ for the similitude groups associated with the following quadratic space (see \cite[\textsection 5]{HK}, \cite{Roberts}).
  For $\aa \in \Qb_{>0}$, let $V = V^{(\aa)}$ be the quadratic space
\begin{equation}
  \label{eq:V}
  V = \left\{
\pmat{a}{\nu}{\nu'}{b}: a, b \in \Qb, \nu \in F
    \right\}
  \end{equation}
  with the quadratic form $Q = Q^{(\aa)}(\lambda) = \aa \cdot \det(\lambda)$.
  We will write $(a, b, \nu)$ to represent $\pmat{a}{\nu}{\nu'}{b} \in V$, and can decompose $V$ into the direct sum of the subspaces
  \begin{equation}
    \label{eq:V01}
    V_0   := \left\{
(0, 0, \nu)
: \nu \in F
\right\} \cong F,~
V_1 := \left\{
  (a, b, 0) : a, b \in \Qb
\right\} \cong \Qb^2. 
\end{equation}
The group $\tH$   acts on $\lambda \in V$ via $\lambda \mapsto a^{-1} \cdot \gamma\lambda (\gamma')^t$ for $(\gamma, a) \in H(\Qb)$, giving rise to $\tH \to
\mathrm{GO}_V$ and $H \cong \mathrm{GSpin}_V  \to \SO_V$ with $H_1 \cong \mathrm{Spin}_V$. 
The GSpin-Shimura variety attached to $V$ is just $X_K$ in \eqref{eq:XK}.

Let $\omega = \otimes_{v \le \infty} \omega_v = \omega_{\psi, V}$ be the Weil representation on $\Sc(V(\Ab)) = \Sc(\hat V) \otimes \Sc(V(\Rb))$, with $\hat V := V(\hat\Qb)$, of $(G \times H)(\Ab)$ associated to $V$.
Following \cite[\S 5]{HK}, we extend it to $R(\Qb)\backslash R(\Ab)$ by setting\footnote{The exponent in (5.1.1) of \cite{HK} should be $-mn/4$ instead of $-mn/2$, in order for (5.1.10) in loc.\ cit.\ to hold.}
\begin{equation}
  \label{eq:omega-R}
  \omega(g, \th)\varphi := \omega(g_1) \cL(\th)\varphi,~
(  \cL(\th)\varphi)(x) := |\nu(\th)|^{-1}\varphi( \th^{-1} \cdot x)
\end{equation}
for $\varphi \in \Sc(V(\Ab))$.
Also note that
\begin{equation}
  \label{eq:center-weil}
\cL(\tilde z)\varphi =   \chi_V(a) \omega(m(a\Nm(z)^{-1})) \varphi 
\end{equation}
  for any $\tilde z = (z, a) \in \Gm(\Ab_F) \times \Gm(\Ab) \cong Z(\tH(\Ab))$ as $\chi_V(\Nm(z)) = 1$.
This allows us to define the theta function 
\begin{equation}
  \label{eq:thetaV}
  \theta(g, \th, \varphi) := \sum_{x \in V(\Qb)} (\omega(g, \th)\varphi)(x)
\end{equation}
on $R(\Qb)\backslash R(\Ab)$, which satisfies the relation
\begin{equation}
  \label{eq:theta-rel}
  \theta(g_1g_2, \th_1\th_2, \varphi) =   \theta(g_1, \th_1, \omega(g_2,\th_2)\varphi)
\end{equation}
for all $(g_i, \th_i) \in R(\Ab)$ \cite[Lemma 5.1.7]{HK}. 

For an even, integral lattice $L \subset V$ with dual $L^\vee$, denote $\Sc_L \subset \Sc(\hat V)$ the finite dimensional subspace spanned by $\cha(\hat L + \mu)$ for $\mu \in \hat L^\vee/\hat L \cong L^\vee/L$ with $\hat L := L \otimes \hat\Zb$.
On $\Sc_L$, we have the following bilinear pairing
\begin{equation}
  \label{eq:bilinear}
  \langle \phi, \varphi \rangle :=
  \sum_{x \in \hat L^\vee/\hat L} \phi(x) \varphi(x). 
\end{equation}

\subsection{Weighted cycle}
\label{subsec:wt-cycle}
In this section, we follow \cite{Kudla97} and recall weighted cycles on Hilbert modular surface $X_K$ attached to the group $H \cong \mathrm{GSpin}_V$ for $V$ from \eqref{eq:V}.

For $x \in V(\Qb)$, we denote its stabilizer by $H_x \subset H$. 
It is isomorphic to $\mathrm{GSpin}(x^\perp)$, which is $B^\times$ for a quaternion algebra $B$ over $\Qb$. 
If $Q(x) > 0$, then $H_x$ gives rise to the natural special cycle $Z(x, h, K)$ on $X_K$ for any $h \in H(\hat\Qb)$ and open compact $K \subset H(\hat\Qb)$. 
Its $\Cb$-points are the image of the map
(see (2.4) in \cite{Kudla97})
\begin{equation}
  \label{eq:Zx}
  H_x(\Qb) \backslash \Db_x \times H_x(\hat\Qb)/K_{x, h}
  \to X_K(\Cb),~
  H_x(\Qb)(z, h_x)K_{x, h} \mapsto
  H(\Qb)(z, h_x h)K_{} 
\end{equation}
with $\Db_x$ the symmetric space for $H_x$ and $K_{x, h} := H_x(\hat\Qb) \cap h K h^{-1}$.
The cycle $Z(x, h, K)$ is defined over $\Qb$ \cite[\textsection 2]{Kudla97}.

Recall we have the subgroup $H_1 \subset H$ defined in \eqref{eq:H1}.
In the decompositions \eqref{eq:X1K} and \eqref{eq:X1K-Gal}, suppose $\delta K \delta^{-1} = K$ for all $\delta \in \Delta$. 
We define the natural special cycle
\begin{equation}
  \label{eq:Z1}
  Z_{1}(x, h, K_1) := Z(x, h, K) \cap X_{1, K_1},  
\end{equation}
which is defined over $E_K$ and the Galois action on $Z_{1}(x, h, K)$'s is the same as in \eqref{eq:E-act}.
In particular
\begin{equation}
  \label{eq:Z1K-Gal}
  Z_1(x, h, K_1)^{\sigma_\delta^{}} \cong Z(x, h, K) \cap X_{1, K, \delta} \cong Z(x, h\delta^{-1}, K) \cap X_{1, K_1} 
\end{equation}
for all $\delta \in \Delta \cong \Cl(\det(K))$. 
Its $\Cb$-points are given as in \eqref{eq:Zx} with $H_x, H$ and $K$ replaced by $H_{1, x}  := H_x \cap H_1, H_1$ and $K_\delta$ respectively.
Using \eqref{eq:X1K} and \eqref{eq:X1K-Gal}, we can write
\begin{equation}
  \label{eq:Z1K}
  Z(x, h, K_1) = \bigsqcup_{\delta \in \Delta} Z_{1}(x, h, K_1)^{\sigma_\delta}
  \subset {X_{K}}_{/E_K}.
\end{equation}

For a $K_1$-invariant Schwartz function $\varphi \in \Sc(\hat V)$ and $x_0 \in V$ with $m := Q(x_0) > 0$, there is a finite set $S(\varphi, x_0) \subset H(\hat\Qb)$ such that \cite[(5.4)]{Kudla97}
\begin{equation}
  \label{eq:Svarphix0}
  \supp \varphi \cap V_m(\hat\Qb)
  =
\bigsqcup_{s \in S(\varphi, x_0)}
  K_1 \cdot s^{-1} \cdot x_0,
\end{equation}
where $V_m := \{x \in V: Q(x) = m\} \subset V$.
The weighted cycle
\begin{equation}
  \label{eq:Z1varphi}
  Z_1(m, \varphi, K_1)
  :=
  \sum_{s \in S(\varphi, x_0)}
  \varphi(s^{-1} \cdot x_0) Z_1(x_0, s , K_1)
\end{equation}
is defined over $E_K$, and independent of the choice of $x_0$ or $S(\varphi, x_0)$.
If $x'_0 = \gamma \cdot x_0$ for $\gamma \in H_1(\Qb)$, then we can choose 
  \begin{equation}
    \label{eq:Svx}
   S(\varphi, x'_0) = \{\gamma \cdot s: s \in S(\varphi, x_0)\}. 
  \end{equation}

\subsection{Theta lifts}
\label{subsec:theta-lifts}
Let $V = V^{(A)}$ be the same as in \eqref{eq:V}, and denote $V(\Ab) = V(\Rb) \oplus \hat V$, where $\hat V := V(\hat\Qb)$. 
We now describe the choice of the archimedean component of the Schwartz function.

In the orthogonal basis  $\{X^+, Y^+, X^-, Y^-\}$ of
  the quadratic space $V(\Rb) = M_2(\Rb)$, where
  \begin{equation}
    \label{eq:ortho-basis}
    Z^+ := X^+ + iY^+ := \pmat{1}{i}{-i}{1},~
    Z^- := X^- + iY^- := \pmat{-1}{i}{i}{1},
  \end{equation}
  For $\epsilon = \pm 1, \delta \in \{+, - \}$, define the following polynomials on $V(\Rb) = M_2(\Rb)$
\begin{equation}
  \label{eq:Pep}
  p^{\delta, \epsilon}(\lambda) :=
(-i) \cdot
  \aa^{-1} (\lambda, X^{-\delta} - i\delta \epsilon Y^{-\delta})
  =
(-i) \cdot  (a - \delta b
+  i \epsilon (\nu + \delta\nu' ))
  \end{equation}
with $\lambda = \smat{a}{\nu}{\nu'}{b} \in V(\Rb)$.\footnote{We will sometimes slightly abuse notation and use $\delta$ to mean $\pm 1$.}
  From these polynomials, we can construct the following Schwartz function  on $V(\Rb)$
  \begin{equation}
    \label{eq:varphi-ep}
    \varphi^{\delta, \epsilon}_{\infty}(\lambda)
    = p^{\delta, \epsilon}(\lambda) e^{-2\pi Q_+(\lambda)} \in \Sc(V(\Rb)),
  \end{equation}
where
    $Q_+(\lambda) :=
Q(\lambda_{Z^+}) - Q(\lambda_{Z^-}) 
=  \frac{\aa}2 \cdot (a^2 + b^2 + \nu^2 + (\nu')^2)$.
  With respect to $K_\infty$, this Schwartz function is right equivariant with weight $(\delta \epsilon, \epsilon)$. 
For $g = g_\tau\in \SL_2(\Rb)$ with  $\tau = u + iv \in \Hb$, we have explicitly
\begin{equation}
  \label{eq:varphi-h-explicit}
\varphi_{\infty, \tau}(\lambda) :=
(\omega_{\infty}(g_\tau)\varphi_\infty)(\lambda)
= v \ebf( Q(\lambda) u) \varphi_\infty^{}(\sqrt{v} \lambda).
\end{equation}

To $\varphi \subset \Sc(\hat V)$ and $\epsilon = \pm 1, \delta \in \{+ , -\}$, we
associate the following theta function for
$(g, h) \in R(\Ab)$ and $\tau \in \Hb$
\begin{equation}
  \label{eq:theta-de}
  \theta^{\delta, \epsilon}(g, h; \varphi) := \theta(g, h; \varphi \otimes \varphi_\infty^{\delta, \epsilon}).
\end{equation}
It is modular of weight $-\delta$ in $g$ and weight $(\delta\epsilon, \epsilon)$ in $h$.
Using computations in the Fock model of $    \omega_\infty$ (see Section 4.1 in \cite{BLY}), we have\footnote{Here, $\stackrel{\cdot}=$ means equality up to multiplication by non-zero constant.}
\begin{equation}
  \label{eq:RL-commute}
  \begin{split}
    R\theta^{1, 1}(g, h, \varphi)
    &\stackrel{\cdot}= R_1 \theta^{-1, 1}(g, h, \varphi)
\stackrel{\cdot}= R_2 \theta^{-1, -1}(g, h, \varphi),\\
    L_1 \theta^{1, 1}(g, h, \varphi)
    &\stackrel{\cdot}= L \theta^{-1, 1}(g, h, \varphi)
  \stackrel{\cdot}= L_2 \theta^{-1, -1}(g, h, \varphi).
  \end{split}
\end{equation}


In $G_F(\Rb) \subset \tH(\Rb)$, we have the following two elements
\begin{equation}
  \label{eq:wi}
  w_1 := \lp \smat{1}{}{}{-1}, \smat{1}{}{}{1} \rp,~
  w_2 := \lp \smat{1}{}{}{1}, \smat{1}{}{}{-1} \rp
\end{equation}
that do not lie in $H(\Rb)$, whereas
$(\dd(-1), w_i) \in R(\Rb)$.
Furthermore, we have
\begin{equation}
  \label{eq:wj-LR}
  R_j(  \rho(w_j) f) =
  \rho(w_j) (L_j f)
\end{equation}
for automorphic form $f$ on $G_F$ and $j = 1, 2$. 
The following observation of their effect on the theta function $\theta^{\delta, \epsilon}$ will be crucial for us later.
\begin{proposition}
  \label{prop:ww'}
  Let $\theta^{\delta, \epsilon}$ be as in \eqref{eq:theta-de}
  and $w_j$ as in \eqref{eq:wi}.
  Then
  \begin{equation}
    \label{eq:theta-ww'}
    \theta^{\delta, \epsilon}(g \dd(-1),  h w_1; \varphi) =
\theta^{-\delta, \epsilon}(g,  h; \varphi),~ 
    \theta^{\delta, \epsilon}(g \dd(-1), h w_2 ; \varphi) =
    \theta^{-\delta, -\epsilon}(g,  h; \varphi)
\end{equation}
for all  $\epsilon = \pm 1, \delta \in \{\pm\}, (g, h) \in R(\Ab)$ and $\varphi \in \Sc(\hat V)$.
\end{proposition}

\begin{proof}
  By definition, we have the following  identity
  \begin{equation}
    \label{eq:varphi-w}
    \begin{split}
(          \omega(\dd(-1), w_1)\varphi^{\delta, \epsilon}_{\infty})( \lambda)
      &=
(\cL(w_1)\varphi^{\delta, \epsilon}_{\infty})(\lambda)        
=\varphi^{\delta, \epsilon}_{\infty}(w_1\cdot \lambda)
        =      \varphi^{-\delta, \epsilon}_{\infty}(\lambda)
    \end{split}
  \end{equation}
  for any $\lambda = \smat{a}{\nu}{\nu'}{b} \in V(\Rb)$, since $w_1\lambda = \smat{a}{\nu}{-\nu'}{-b}$ and
  $$
  p^{\delta, \epsilon}(\smat{a}{\nu}{-\nu'}{-b}) 
  =
(-i) \cdot
  (a + \delta b
    +  i \epsilon (\nu - \delta \nu' ))
    =   p^{-\delta, \epsilon}(\smat{a}{\nu}{\nu'}{b}) 
$$
by \eqref{eq:Pep}.
Similarly, $\omega(\dd(-1), w_2)\varphi^{\delta, \epsilon}_{\infty} = \varphi^{-\delta, -\epsilon}_{\infty}$. 
As a result of \eqref{eq:theta-rel}, we obtain  \eqref{eq:theta-ww'}.
\end{proof}

Let $f_0 \in S_1(N, \chi_0)$ be a holomorphic elliptic modular form with adelization  $\fs_0: G(\Ab) \to \Cb$.
For $\varphi \in \Sc(\hat V)$, 
we define the theta lift $\Phi^{}(h; f_0, \varphi)$ of $f_0$ by
\begin{equation}
  \label{eq:Phi}
  \Phi^{}( h; f_0, \varphi)
  =       \Phi^{}( h; f^\#_0, \varphi)
      :=
      \int_{[G]}  \fs_0(g \dd(\nu(h))) \theta^{1, 1}(g \dd(\nu(h)),  h; \varphi)  dg.
    \end{equation}
    for $h \in H(\Ab)$. 
    It is left $H(\Qb)$-invariant since
    \begin{align*}
          \Phi^{}(\gamma h; f_0, \varphi)
 			&  = \int_{[G]}  \fs_0(g \dd(\nu(\gamma h))) \theta^{1, 1}(g \dd(\nu(\gamma h)),  \gamma h; \varphi) dg\\
 			&  = \int_{[G]}  \fs_0(g' \dd(\nu( h))) \theta^{1, 1}(g' \dd(\nu(h)),   h; \varphi) dg \\
 			&  = \Phi^{}( h; f_0, \varphi)
    \end{align*}
for $\gamma \in H(\Qb)$ with $g' = \gamma^{-1} g \gamma \in G(\Ab)$.
It is easy to verify that
\begin{equation}
  \label{eq:rho-comm}
  \Phi^{}( h; \rho(g_1) f_0, \varphi)
  =
   \Phi^{}( h; f_0, \omega(g^{-1}) \varphi). 
 \end{equation}
 for any $g \in G(\hat\Qb)$. 

  Similarly given a Hilbert cusp form $f$ of weight $(1, 1)$, 
  we define its theta lift
\begin{equation}
  \label{eq:I}
  I(g; f, \varphi) :=
  \int_{[H_1]}  f( h_1 h)
  \overline{\theta^{1, 1}(g,  h_1h; \varphi)} dh_1
\end{equation}
for $g \in \tG(\Ab)$ with $h \in \tH(\Ab)$ satisfying $\nu(h) = g$.
It is easy to see that
  \begin{equation}
    \label{eq:lift-comm}
    \int_{[G]} \overline{I(g; f, \varphi) }   f_0^\#(g) dg
    =
    \int_{[H_1]} \Phi^{}(h; f_0, \varphi) f(h) dh.
  \end{equation}

  \begin{lemma}
    \label{lemma:lift-hol}
If $f_0$, resp.\ $f$, is holomorphic, then its theta lift $\Phi$, resp.\ $I$, is holomorphic. 
  \end{lemma}

  \begin{proof}
    This follows directly from \eqref{eq:RL-commute} and Stokes' theorem.
  \end{proof}

\section{Hecke equivariance of Doi--Naganuma lifts}
\label{subsec:Hecke}
In this section, we define Hecke operators in $\tG(\Ab)$ and $\tH(\Ab)$, and show that the theta lift from Section \ref{subsec:theta-lifts} intertwines the actions in a suitable way. Such result are known in the literature as generalizations of the Eichler commutation relation \cite{Eichler}. For Hecke operators on $G \times H$, the relation is worked out by Rallis in the representation language \cite{Rallis82}. For our purpose, it is necessary to work with Hecke operators on similitude groups, which does not seem to be in the literature. Therefore, we give the detailed proof here, in a rather classical way that transfers all the actions to the Schwartz function.

In the notation before, let $f$ be an automorphic form of parallel weight $k$, level $N$ on $\tG(\Ab) = \GL_2(\Ab_F)$, and $\pf$ a prime of $F$ such that it is co-prime to $DN$ and $f$ is right $\tH(\Zb_p)$-invariant, with $p$ the rational prime below $\pf$.
For  $K_\pf$-biinvariant $\phi \in \Cc(\tH(\Qb_p))$, where  $K_\pf := \GL_2(\Oc_\pf)$ with $\Oc_\pf$ the valuation ring of $F_\pf$, 
the associated Hecke operator $\Tc_\phi$ is defined by
\begin{equation}
  \label{eq:Tphi}
  (  \Tc_\phi f)(h)
  := \int_{\tH(\Qb_p)} f(h \tilde h^{-1}) \phi(\tilde h) d \tilde h. 
\end{equation}
If $\phi = \cha(K_\pf \varpi \dd(1/\varpi) K_\pf)$ with $\varpi \in \Oc_{\pf}$ a uniformizer, we denote $\Tc_\phi$ by $\Tc_\pf$. 
It is explicitly given by
\begin{equation}
  \label{eq:Tcp}
  (  \Tc_\pf f)(h) =
  \sum_{\beta \in \dd(\varpi)K_\pf\dd(\varpi^{-1})) \cap K_\pf\backslash K_\pf}
  f(h (\varpi\dd(1/\varpi) \beta)^{-1}).
\end{equation}
In addition, we also have the diamond operator
\begin{equation}
  \label{eq:diamond}
(  \langle \pf \rangle f)(h) := f( h \varpi^{-1}). 
\end{equation}
Suppose $F$ has narrow class number 1 and $f$ has weight $(k, k)$. Then it is easy to verify that 
\begin{equation}
  \label{eq:Hecke-compare}
(\Tc_\pf f)^\flat
= \Nm(\pf)^{1-k/2}   T_\lambda f^\flat ,~
(\langle \pf \rangle f)^\flat =
\langle \lambda \rangle f^\flat,
\end{equation}
where $\pf = (\lambda)$ and $T_\lambda, \langle \lambda\rangle$ are the classical Hecke and diamond operators.
It turns out that the action of $\Tc_\pf$ on the theta lift $\Phi^{\delta, \epsilon}$ defined in \eqref{eq:Phi} is compatible with the Hecke operator on the input.

\begin{proposition}
  \label{prop:Hecke-comm}
  Let $f_0 \in S_1(N, \chi_0)$ and
  $\pf$ be a prime of $F$ above rational $p$ co-prime to $\aa DN$.
  For $\varphi = \varphi^{(p)} \otimes_{} \phi \in \Sc(\hat V)$ with $\varphi^{(p)} \in \Sc(\otimes_{p \neq \ell < \infty} V_\ell)$ and $\phi = \cha(L)$ with $L := M_2(\Oc_{ \pf}) \cap V_p$, we have
  \begin{equation}
    \label{eq:Hecke-comm}
    \begin{split}
          \Tc_\pf \Phi^{\delta, \epsilon}(h; f_0, \varphi)
&    =
    \begin{cases}
      \Phi^{\delta, \epsilon}(h; \Tc_pf_0, \varphi) & \text{ if } t = 1,\\
      \Phi^{\delta, \epsilon}(h; (\Tc^2_p - 2 p \langle p^{} \rangle )f_0, \varphi) & \text{ if } t = 2,
    \end{cases}\\
          \langle \pf \rangle \Phi^{\delta, \epsilon}(h; f_0, \varphi)
&    =
 \Phi^{\delta, \epsilon}(h; \langle q\rangle f_0, \varphi),
    \end{split}
  \end{equation}
  where $q = |\Oc_{ \pf}/\varpi| = p^t$. 
\end{proposition}

\begin{remark}
  \label{rmk:center}
  The second equation in \eqref{eq:Hecke-comm} is a special case of the relation
  \begin{equation}
    \label{eq:center-act}
    \Phi^{\delta, \epsilon}(h z; f_0, \varphi) =
        \Phi^{\delta, \epsilon}(h ; \rho(\Nm(z))f_0, \varphi) 
  \end{equation}
  for any $z \in Z_F(\Ab)$, which follows from \eqref{eq:center-weil}.
\end{remark}
\begin{proof}
For $r \in \Zb/q\Zb$,   we define the following Schwartz functions in $\Sc_{L}$
\begin{equation}
  \label{eq:Lr}
  \begin{split}
    \phi_{r} &:= \cha(q^{-1} L_{r}),~  
  L_{r} := \{\lambda \in L - p L: \det(\lambda) \equiv r \bmod q\} \subset L. 
  \end{split}
\end{equation}
The coset in \eqref{eq:Tcp} has size $q + 1$ and is given by
  \begin{equation}
    \label{eq:B}
      \dd(\varpi)K_\pf\dd(\varpi^{-1})) \cap K_\pf\backslash K_\pf
      \cong \Gamma_0(\pf) \backslash \SL_2(\Oc)
  \cong \{ n^-(j): j \in \Oc_\pf/\pf \} \cup \{w\} =:B ,
  \end{equation}
  where
  $w = \smat{}{-1}{1}{},   n^-(j)
  =   \smat{1}{0}{j}{1}$.
  We set $\tphi :=  \cL(\varpi^{-1}\dd(\varpi)) \phi
  $ and need to evaluate 
    \begin{equation}
      \label{eq:tvp-ave}
 \sum_{\beta \in B} \cL(\beta^{-1} \varpi^{-1} \dd(\varpi)) \phi
  =
  \sum_{\beta \in B}
\omega(\beta^{-1})\tphi
    \end{equation}
This is because
$$
    \Phi^{\delta, \epsilon}(h (\varpi\dd(1/\varpi)\beta)^{-1}; f_0, \phi)
    =
\Phi^{\delta, \epsilon}(h; \rho(\dd(1/p))f_0,
\cL(\beta^{-1}\varpi^{-1}\dd(\varpi))      \phi)
$$
  
  We first treat the case $q = p$, i.e.\ $p = \pf\pf'$ is split.
  Then $F_\pf \cong \Qb_p$ identifies $V_p \cong M_2(\Qb_p), K_\pf \cong \GL_2(\Zb_p)$, under which $\dd(\varpi) \in H(\Qb_p)$ and $\phi$ become $(\dd(p), 1) \in \GL_2(\Qb_p)^2$ and $\cha(M_2(\Zb_p))$ respectively.
  Then $\tphi =
p^{-1}    \cha\lp \smat{p^{-1}\Zb_p}{p^{-1}\Zb_p}{\Zb_p}{\Zb_p}  \rp$ and we claim that
    \begin{equation}
      \label{eq:tvp-ave-split}
  \sum_{\beta \in B} \cL(\beta^{-1}\varpi^{-1} \dd(\varpi)) \phi
  =
  \sum_{\beta \in B}
\omega(\beta^{-1})\tphi
  = (\phi_{0} + (p+1) \phi_{})/p.
    \end{equation}
      Indeed the left hand side has support on $L$, is locally constant with respect to translation by $L$. Therefore it is contained in $\Sc_{L}$. 
Furthermore, it is invariant with respect to $K_{1, p} := \SL_2(\Zb_p)^2 \subset H_1(\Qb_p)$. 
The subspace of such function is spanned by $\phi_{r}$ and $\phi_{}$, which are orthogonal with respect to the pairing in \eqref{eq:bilinear}.
So we can write
  \begin{align*}
  \sum_{\beta \in B}
&\omega(\beta^{-1})\tphi
  =
  \sum_{r \in \Zb/p\Zb}
\frac{  \left\langle
  \sum_{\beta \in B}\omega(\beta^{-1}) \tphi, \phi_{r} \right\rangle}{  \langle \phi_{r}, \phi_{r}\rangle} \phi_{r}
+
  \left\langle     \sum_{\beta \in B} \omega(\beta^{-1}) \tphi, \phi_{} \right\rangle \phi_{}\\
&= |B| \cdot
  \lp
  \sum_{r \in (\Zb/p\Zb)^\times}
\frac{  \langle
\tphi, \phi_{r} \rangle}{ (p^2-1)p} \phi_{r}
   +
  \frac{  \langle
\tphi, \phi_{0} \rangle}{ (p^2-1)(p+1)} \phi_{0}
+\langle \tphi, \phi_{} \rangle \phi_{}
  \rp
      \end{align*}
      It is easy to see that $\langle
      \tphi, \phi_{} \rangle = 1/p$ and
      $$
\langle
\tphi, \phi_{r} \rangle
=
\begin{cases}
\tfrac{  p^2 - 1}p & \text{if } r = 0,\\
  0 & \text{otherwise.}
\end{cases}
      $$
  Applying $|B| = p+1$ then gives us \eqref{eq:tvp-ave-split}, and we can write
  \begin{equation}
    \label{eq:Tcp-split}
    (  \Tc_\pf \Phi^{\delta, \epsilon})(h; f_0, \varphi)
    =
      \Phi^{\delta, \epsilon}(h; \rho(\dd(p^{-1}))f_0,       \varphi^{(p)} \otimes
(\phi_{0} + (p+1)\phi_{})/p).
  \end{equation}
  On the other hand, we have
  \begin{equation}
    \label{eq:Tp}
    \Tc_p      \fs_0
    =
\sum_{\beta \in B}        \rho(\beta^{-1}p^{-1}\dd(p))\fs_0.
  \end{equation}
  Applying Equation \eqref{eq:rho-comm}, we have
  \begin{align*}
    \Phi^{\delta, \epsilon}
    &(h; \Tc_p f_0, \varphi)
    =
      \sum_{\beta \in B}
      \Phi^{\delta, \epsilon}(h;       \rho(\beta^{-1}p^{-1}\dd(p))f_0, \varphi)
      =
      \sum_{\beta \in B}
      \Phi^{\delta, \epsilon}(h;       \rho(p^{-1}\dd(p))f_0, \omega(\beta)\varphi)      \\
    &= (p+1)       \Phi^{\delta, \epsilon}(h;       \rho(p^{-1}\dd(p))f_0,\varphi) .
  \end{align*}
  Now since $\rho(\kappa) f_0 = f_0$ for any $\kappa \in \SL_2(\Zb_p) \subset G(\Ab)$, we have
  \begin{align*}
      \Phi^{\delta, \epsilon}(h;       \rho(p^{-1}\dd(p))f_0,\varphi)
&  =
  \Phi^{\delta, \epsilon}(h;       \rho(w^{-1} \dd(p^{-1}) w)f_0,\varphi)
  =
      \Phi^{\delta, \epsilon}(h;       \rho(\dd(p^{-1}))f_0,\varphi)
  \end{align*}
  Similarly for any $j \in \Qb_p$, we have
  \begin{align*}
      \Phi^{\delta, \epsilon}&(h;       \rho(p^{-1}\dd(p))f_0,\varphi)
  =
  \Phi^{\delta, \epsilon}(h;       \rho(p^{-1}\dd(p)n(j))f_0,\varphi)  \\
&  =
  \Phi^{\delta, \epsilon}(h;       \rho(n(j/p)m(p^{-1})\dd(p^{-1}))f_0,\varphi)  
  =
  \Phi^{\delta, \epsilon}(h;       \rho(\dd(p^{-1}))f_0,\omega(m(p) n(-j/p))\varphi)  .
  \end{align*}
  Averaging this over $j \in \Zb/p\Zb$ gives us
  $$
  \Phi^{\delta, \epsilon}(h;       \rho(p^{-1}\dd(p))f_0,\varphi)
  =
      \Phi^{\delta, \epsilon}(h;       \rho(\dd(p^{-1}))f_0,\varphi^{(p)} \otimes (\phi_{0} + \phi)/p^2)  .
  $$
  Putting these together give us \eqref{eq:Hecke-comm}.

  For the case $\pf = p$ is inert, we have $\varpi = p, q = p^2$.
  In addition to the Schwartz function in \eqref{eq:Lr}, we define
  \begin{equation}
    \label{eq:phi-inert}
    \begin{split}
      \phi'
      &:= \cha(\frac1p L) = \phi + \sum_{r \in \Zb/p\Zb} \phi_{p, r}, \\
    \phi_{p, r} &:= \cha(\frac1p L_{p, r}),~ L_{p, r} :=
    \{ \lambda: \lambda \in L - pL, \det(\lambda) \equiv r \bmod p\} \subset L - p L. 
    \end{split}
  \end{equation}
  for $r \in \Zb/p\Zb$.
  Then $\tphi = q^{-1} \cha(\smat{p^{-2}\Zb_p}{p^{-1}\Oc}{p^{-1}\Oc}{\Zb_p} \cap V_p)$ and the analogue of \eqref{eq:tvp-ave-split} is
  \begin{equation}
    \label{eq:tvp-ave-inert}
    \begin{split}
q        \sum_{\beta \in B}
\omega(\beta^{-1})\tphi
      &  = \phi_{0} + (q+1) \phi_{} + \phi_{p, 0} + (p+1)\sum_{r \in (\Zb/p\Zb)^\times} \phi_{p, r}\\
      &= \phi_{0} + (q-p)\phi_{}. + (p+1)\phi' - p\phi_{p, 0} 
    \end{split}
  \end{equation}
  So we have
  \begin{equation}
    \label{eq:Tcp-inert}
        (  \Tc_\pf \Phi^{\delta, \epsilon})(h; f_0, \varphi)
    =
q^{-1}      \Phi^{\delta, \epsilon}(h; \rho(\dd(q^{-1}))f_0,       \varphi^{(p)} \otimes
(\phi_{0}  + (q-p)\phi_{}+ (p+1)\phi' - p\phi_{p, 0})).
\end{equation}

On the other hand, we have
\begin{align*}
  \Phi^{\delta, \epsilon}(h;       \Tc_p^2 f_0,\varphi)
  &= (p+1)   \Phi^{\delta, \epsilon}(h;       \rho(p^{-1}\dd(p)) (\Tc_pf_0),\varphi)\\
&    = (p+1)
  \Phi^{\delta, \epsilon}(h;       \rho(\dd(q^{-1})) f_0, \varphi^{(p)} \otimes (p \phi - p^{-2} \phi' ) ),\\
  \Phi^{\delta, \epsilon}(h;          \langle p\rangle f_0,\varphi)
  &=  - p^{-2} \Phi^{\delta, \epsilon}(h;       \rho(\dd(q^{-1})) f_0, \varphi^{(p)}\otimes \phi' ).
\end{align*}
Using the action $\omega(m(p)) \phi = - p^{-2} \phi'$
and 
the right, resp.\ left, $\SL_2(\Zb_p)$-invariance of $f_0$, resp.\ $\phi$, one can deduce that
\begin{align*}
  \Phi^{\delta, \epsilon}(h;          \rho(\dd(q^{-1})) f_0,\varphi^{(p)} \otimes \phi_{p, 0})
  &=  \Phi^{\delta, \epsilon}(h;       \rho(\dd(q^{-1})) f_0, \varphi^{(p)} \otimes ( (-p^2+p-1 )\phi + \tfrac{1}{p}  \phi' )\\  
  \Phi^{\delta, \epsilon}(h;          \rho(\dd(q^{-1})) f_0,\varphi^{(p)} \otimes \phi_{0})
  &=  \Phi^{\delta, \epsilon}(h;       \rho(\dd(q^{-1})) f_0, \varphi^{(p)} \otimes  (p^4 \phi - \phi')).
\end{align*}
Putting these together into \eqref{eq:Tcp-inert} gives \eqref{eq:Hecke-comm} for $e = 2$. 
\end{proof}

Using Equation \eqref{eq:lift-comm} and the adjointness of Hecke operator with respect to the Petersson inner product, we have the following corollary.

\begin{corollary}
  \label{cor:Hecke-comm2}
  Given a Hilbert cusp form $f$ of level $\Nc$, weight $(1, 1)$ 
  and prime $\pf$ with norm $q = p^t$ co-prime to $A \Nm(\Nc) D$, we have
  \begin{equation}
    \label{eq:Hecke-comm2}
    (    \Tc_q - 2 p \langle p^{-1} \rangle) I(g_0, f, \varphi) =
 I(g_0, \Tc_\pf f, \varphi).
  \end{equation}
\end{corollary}

\section{Fourier expansions of Doi-Naganuma lifts}

In this section, we will describe the data for our theta lift. 
Throughout, $V = V^{(D)}$ is the same as in \eqref{eq:V}, with odd $D > 0$, the discriminant of a real quadratic field $F$. We fix a normalized newform $f_0 \in S_1(N, \chi_0)$ with $N$ square-free and co-prime to $2D$.
Denote by $f: [G_F] \to \Cb$ its base-change, normalized with $c_{1, \df}(f) = 1$.
It is right $K(N, \df)$-equivariant with character $\chi := \chi_0 \circ \Nm_{F/\Qb}$, where
\begin{equation}
  \label{eq:KN}
  K(N, \df) := \{\gamma \in
  \GL_2(\hat\Oc): \gamma \equiv \smat{*}{*}{0}{*} \bmod{N},~ \det(\gamma) \equiv 1 \bmod \df\}
\end{equation}
is an open compact subgroup of $\GL_2(\hat\Oc)$ and $\chi(\smat{\alpha}***) := \chi(\alpha)$.

\subsection{Schwartz function}
\label{subsec:schwartz}
We start with the Schwartz function. 
 For a character $\chi$ of $(\Zb/p\Zb)^\times$, 
let  $\gf_\chi := \sum_{j \in (\Zb/p\Zb)^\times} \chi(j) \ebf(j/p)$ be its Gauss sum. 
Define $\varphi_p \in \Sc(V_p)$ by
\begin{equation}
  \label{eq:varphip}
  \varphi_p =
  \begin{cases}
    \cha (\Zb_p^2 \oplus \Oc_p) & \text{ if } p \nmid DN,\\
\frac1{\gf_{\overline{\chi_{0, p}}}}
    \sum\limits_{j \in (\Zb/p\Zb)^\times} \overline{\chi_{0,p}(j)}
    \lp    \cha(\frac{j}p + \Zb_p) - \frac1p \cha(p^{-1}\Zb_p) \rp
    \otimes \cha(p\Zb_p \oplus \Oc_p)
                               & \text{ if } p \mid N,\\
    \frac1{\gf_{\chi_{F, p}}p(p-1)}    \sum\limits_{h \in \GL_2(\Oc_p/\df_p)} \chi_{F, p}(\det(h)) \cha(h \cdot
    ((\frac1p, 0, 0) + \Zb_p^2 \oplus \df_p^{-1}))
                               & \text{ if } p \mid D.
  \end{cases}
\end{equation}
In the case $p \mid D$, denote $L := \Zb_p^2 \oplus \df_p^{-1}$ with the quadratic form $Q(a, b, \nu) = p(ab-\Nm(\nu))$. 
For $\mu \in L^\vee/L$, the intersection $\supp(\varphi_p) \cap L + \mu$ is non-empty if and only if $\mu$ is an isotropic vector, i.e.\ $Q(\mu) = 0 \in \Qb/\Zb$ and $\mu \neq 0$. 
Since $D$ is odd, we have $\GL_2(\Oc_p/\df_p) \cong \GL_2(\Zb/p\Zb)$
acting on $L^\vee/L \cong (\frac1p\Zb/\Zb)^3 \cong \mathrm{Sym}_2(\Zb/p\Zb)$ by conjugation and scaling by determinant.
The set of isotropic vectors in  are given by
$$
\{\tfrac{d}p(1, 0, 0): d \in (\Zb/p\Zb)^\times\}  \cup
\{\tfrac{d}p(j^2, 1, j): j \in \Zb/p\Zb, d \in (\Zb/p\Zb)^\times\},
$$
since it is acted on by $\GL_2(\Zb/p\Zb)$ transitively and the stabilizer of $(\frac1p, 0, 0)$ is generated by $\smat1101$ and the center. 
So we can explicitly write
\begin{equation}
  \label{eq:varphipD}
  \varphi_p =
  \sum_{d \in (\Zb/p\Zb)^\times}
  \chi_{F, p}(d)
\lp  \cha(\tfrac{d}p(1, 0, 0) + L) + \sum_{j \in \Zb/p\Zb}   \cha(\tfrac{d}p(j^2, 1, j) + L)\rp
\end{equation}
when $p \mid D$.

It is now easy to check that $\varphi_p$ satisfies the following result.
\begin{lemma}
\label{lemma:varphi-cond}
  The Schwartz function $\varphi = \otimes_{p < \infty} \varphi_p$ defined in \eqref{eq:varphip} satisfies 
  \begin{equation}
  \label{eq:varphi-cond}
  \omega_f(\kappa, h) \varphi = \chi_0(\kappa)^{} {\chi(h^{})}
  \varphi,~
\end{equation}
for all $\kappa \in K_0(N) \subset \SL_2(\hat\Zb) \subset G(\hat\Qb)$ and $h \in K_{}(N, \df) \subset \GL_2(\hat\Oc)$.

\end{lemma}
\begin{proof}
  We check this claim locally.
  When $p\nmid DN$,  the Schwartz function $\varphi_p$ is spherical and the claim is clear.
  When $p \mid D$, we have $\omega_p(\kappa) \varphi_p = \varphi_p$ for all $\kappa \in \SL_2(\Zb_p)$ since $\varphi_p$ is an invariant vector \cite[\textsection 5.1]{LZ-DN}.
  When $p \mid N$ and $\kappa = m(a) \in K_0(N) \cap \SL_2(\Zb_p)$, we have
  $$
(\omega_p(\kappa) \varphi_p)(x) = \varphi_p(ax) = \chi_{0, p}(a) \varphi_p(x). 
  $$
  This proves the claim for $\kappa \in K_0(N)$.
  The claim for $h \in K(N, \df)$ follows easily by substitution. This completes the proof. 
\end{proof}

Now let $\ell = (1, 0, 0), \ell' = (0, {1/D}, 0) \in V \subset V_p$ be isotropic vectors, which satisfy $(\ell, \ell') = D$.
They give the following partial Fourier transform
\begin{equation}
  \label{eq:Fc}
  \Fc(\varphi)(\eta, \nu)
  = \int_\Ab \varphi(a \ell + \eta_1 \ell' + (0, 0, \nu))\psi(a\eta_2 ) da
\end{equation}
for $\varphi \in \Sc(V(\Ab))$ and $\eta = (\eta_1, \eta_2) \in \Ab^2$.
This linearizes the symplectic part of the Weil representation on $V_1 \subset V$, i.e.\
\begin{equation}
  \label{eq:linearize}
  \Fc(\omega_1(g)\varphi)(\eta, \nu)
  =   \Fc(\varphi)(\eta \cdot g, \nu)
\end{equation}
for all $g \in G(\Ab)$. Furthermore, Poisson summation gives us
\begin{equation}
  \label{eq:Poisson}
  \theta(g, h, \varphi) = \sum_{x \in V} (\omega(g, h)\varphi)(x)
  = \sum_{\eta \in \Qb^2,~ \nu \in V_0} \Fc(\omega(g, h)\varphi)(\eta, \nu).
\end{equation}
If $\varphi = \otimes_{p \le \infty}\varphi_p$, then $\Fc(\varphi) = \prod_{p \le \infty} \Fc_p(\varphi_p)$, where $\Fc_p$ is the Fourier transform on $\Sc(V_p)$ defined similarly as in \eqref{eq:Fc}. 

We record the partial Fourier transforms of the Schwartz functions $\varphi_p$ defined in \eqref{eq:varphip}. This will be used
in Section \ref{subsec:unfoldI} when we compute the Fourier expansion of the theta lift.

\begin{proposition}
  \label{prop:PFT}
  Let $\varphi_p$ be defined as in \eqref{eq:varphip}. Then its partial Fourier transform as defined in \eqref{eq:Fc} is given by
  \begin{equation}
    \label{eq:PFT-varphip}
    \Fc_p(\varphi_p)      (\eta, \nu)
    =
    \begin{cases}
      \cha(\Zb_p^2 \oplus \Oc_p)      (\eta, \nu) & \text{if } p \nmid DN,\\
      \chi_{0, p}(-\eta_2)
      \cha(p\Zb_p \oplus \Zb_p^\times \oplus \Oc_p)       (\eta, \nu)& \text{if }p \mid N,\\
 \chi_{F, p}(-\eta_2)      \cha(p \Zb_p \oplus \Zb_p^\times \oplus \df^{-1}_p)      (\eta, \nu)&\\
      +
    \frac1{\gf_{\chi_{F, p}}}
      \displaystyle\sum_{\substack{d \in (\Zb/p\Zb)^\times\\ j \in \Zb/p\Zb}} \chi_{F, p}(d) \ebf(-\tfrac{d j^2\eta_2}p)\\
      \qquad \times \cha((p\Zb_p + d) \oplus \Zb_p \oplus (\df^{-1}_p + \tfrac{dj}{p}))
      (\eta, \nu)
                                 & \text{if } p \mid D.
    \end{cases}
  \end{equation}
\end{proposition}
\begin{remark}
  \label{rmk:PFT0}
  When $\eta = (0, \eta_2)$, we have
  \begin{equation}
    \label{eq:PFT0}
    \Fc_p(\varphi_p)      (\eta, \nu)
    =
    \begin{cases}
      \cha(\Zb_p \oplus \Oc_p)      (\eta_2, \nu) & \text{if } p \nmid DN,\\
      \chi_{0, p}(-\eta_2)
      \cha(\Zb_p^\times \oplus \Oc_p)       (\eta_2, \nu)& \text{if }p \mid N,\\
      \chi_{F, p}(-\eta_2)      \cha(\Zb_p^\times \oplus \df^{-1}_p)      (\eta_2, \nu) & \text{if } p \mid D. 
\end{cases}
\end{equation}
This implies that $\Fc(\varphi)(\eta r, \nu) = \chi_0(r) \chi_F(r) \Fc(\varphi)(\eta, \nu)$ for any $r \in \hat\Zb^\times$. 
Also for $p \mid N$, $\Zb_p \times \{0\}$ does not intersect the support of $\varphi_p$.
\end{remark}

\begin{proof}
  When $p \nmid D$, i.e.\ in the first 2 cases, we can write $\varphi_p = \varphi_{1, p} \otimes \cha(\Oc_p)$, where $\varphi_{1, p} \in \Sc(V_{1, p})$. Then $\Fc_p(\varphi_p) = \Fc_p(\varphi_{1, p}) \otimes \cha(\Oc_p)$ with
  $    \Fc_p(\varphi_{1, p}) = \int_{\Qb_p}\varphi_{1, p}(a \ell + \eta_1 \ell')\psi_p(a \eta_2) da$.
  When $p \nmid DN$, we have $\varphi_{1, p} = \cha(\Zb_p^2) = \Fc_p(\varphi_{1, p})$.
  When $p \mid N$, we have
  \begin{align*}
\gf_{\overline{\chi_{0, p}}}    \Fc_p(\varphi_{1, p})(\eta)
    &=
      \cha(p\Zb_p)(\eta_1)
      \sum_{j \in (\Zb/p\Zb)^\times} \overline{\chi_{0,p}(j)}
\lp      \int_{\frac{j}p + \Zb_p}\psi_p(a \eta_2) da
      -\frac1p \int_{p^{-1}\Zb_p}\psi_p(a \eta_2) da \rp\\
    &=
      \cha(p\Zb_p)(\eta_1)
      \sum_{j \in (\Zb/p\Zb)^\times} \overline{\chi_{0,p}(j)}
\lp \ebf(-\frac{\eta_2 j}p)\cha(\Zb_p)
      - \cha(p\Zb_p)\rp(\eta_2)\\
    &=
            \cha(p\Zb_p \oplus \Zb_p^\times)(\eta)
\chi_{0, p}(-\eta_2)\gf_{\overline{\chi_{0, p}}}. 
  \end{align*}
  Here we have used $\ebf(-\frac{\eta_2 j}p)\cha(\Zb_p)(\eta_2)
      - \cha(p\Zb_p)(\eta_2) = \ebf(-\frac{\eta_2 j}p)\cha(\Zb_p^\times)(\eta_2)$. 
  When $p \mid D$, we use \eqref{eq:varphipD} to obtain
  \begin{align*}
\gf_{\chi_{F, p}}    \Fc_p(\varphi_p)(\eta, \nu)
    &=
  \sum_{d \in (\Zb/p\Zb)^\times}
  \chi_{F, p}(d)
\lp
      \begin{aligned}[c]
&      \cha(p\Zb_p)(\eta_1)
\cha(\df_p^{-1})(\nu)
  \int_{\tfrac{d}p + \Zb_p} \psi_p(a \eta_2) da\\
&
+      \cha(d+p\Zb_p)(\eta_1)
  \sum_{j \in \Zb/p\Zb}
\cha(\tfrac{dj}p + \df_p^{-1})(\nu)
      \int_{\tfrac{dj^2}p + \Zb_p} \psi_p(a \eta_2) da 
      \end{aligned}
      \rp\\
    &= 
  \sum_{d \in (\Zb/p\Zb)^\times}
  \chi_{F, p}(d)
      \ebf(-\tfrac{d\eta_2}p)
      \cha(p\Zb_p\oplus \Zb_p \oplus  \df_p^{-1})(\eta, \nu)\\
    & \quad + 
  \sum_{d, j \in (\Zb/p\Zb)^\times}
  \chi_{F, p}(d)
            \ebf(-\tfrac{dj^2\eta_2}p)
      \cha((d + p\Zb_p) \oplus\Zb_p \oplus (\tfrac{dj}p + \df_p^{-1}))(\eta, \nu).
  \end{align*}
  This finishes the proof.
\end{proof}

\subsection{Lifting identification.}
The goal of this section is to identify the theta lifts of $f$, resp.\ $f_0$, with $f_0$, resp.\ $f$.
\begin{lemma}
  \label{lemma:wt1-id}
  For $i = 1, 2$, let $f_i \in S_1(N_i, \chi_i)$ be normalized newforms with Fourier coefficients $\{a_i(n): n \ge 1\}$ and associated Galois representations $\varrho_i: G := \Gal(L/\Qb) \to \GL_2(\Cb)$ for some number field $L$.
  If every conjugacy class of $G$ contains $\mathrm{Frob}_p$ for a prime $p$ such that $a_1(p) = a_2(p), \chi_1(p) = \chi_2(p)$, then $f_1 = f_2$. 
\end{lemma}

\begin{proof}
Under the assumption, we have $\varrho_1 \cong \varrho_2$, which implies $f_1 = f_2$. 
\end{proof}
Now let $f_0 \in S_1(N, \chi_0)$ be a normalized newform. 
For a real quadratic field $F=\Qb(\sqrt{D})$ with fundamental discriminant $D$, denote $f: [H] \to \Cb$ the normalized base-change of $f_0$ to $F$ of weight $(1, 1)$ and character $\chi(\alpha) = \chi_0(\Nm(\alpha))$ for $\alpha \in \Ab_F^\times$. 
Let $\varrho_f$ and $\Q(f)$, resp.\ $\varrho_{f_0}$ and $\Q(f_0)$, denote the Artin representation associated with $f$, resp.\ $f_0$, and number field generated by the Fourier coefficients of $f$, resp.\ $f_0$. Then $\varrho_f$ is just the restriction of $\varrho_{f_0}$ to $\Gal(\overline{F}/F) \subset \Gal(\overline{\Q}/\Q)$, from which it is clear that $\Q(f) \subset \Q(f_0)$. In general, the latter is at most a quadratic extension of the former. In the following case, they turn out to be the same.

\begin{lemma}
  \label{lemma:Qf}
  Suppose $D$ and $N$ are co-prime. Then $\Q(f) = \Q(f_0)$. 
\end{lemma}

\begin{proof}
  Let $L/\Q$ be the fixed field of the kernel of $\varrho_{f_0}$. Then $L$ and $F$ are disjoint since they are ramified at primes dividing $N$ and $D$ respectively.
  By the Chebotarev density theorem, every conjugacy class in $\Gal(L/\Q)$ contains $\mathrm{Frob}_p$ with $p$ splits in $F$. In that case, the trace of $\varrho_f(\mathrm{Frob}_{p}) =  \varrho_{f_0}(\mathrm{Frob}_{p})$ is contained in $\Q(f)$. 
  Since $\Q(f_0)$ is generated by $\tr(\varrho_{f_0}(C))$ over conjugacy classes  $C \subset \Gal(L/\Q)$, we  obtain $\Q(f_0)  \subset \Q(f)$. 
\end{proof}

For $\varphi \in \Sc(\hat V)$, we have the following result.
\begin{proposition}
  \label{prop:rr0}
  There exists constants $r = r(\varphi), r_0 = r_0(\varphi)$ in $\Cb$ such that
  \begin{equation}
    \label{eq:lift-eq}
    \Phi^{1, 1}(h; f_0, \varphi) = r\cdot f(h),~
    I(g; f, \varphi) = r_0\cdot f_0^\#(g).
  \end{equation}
  Furthermore, we have
  \begin{equation}
    \label{eq:rr0}
    \overline{r_0} \| f_0\|_{1, \mathrm{Pet}}^2 =     r \| \rho(\xi) f\|_{1, \mathrm{Pet}}^2,
  \end{equation}
  for all $\xi \in H(\Ab)$. 
\end{proposition}

\begin{remark}
  \label{rmk:rr0}
  Since $f_0, f$ are normalized, the constants $r_0, r$ are the first Fourier coefficients of $I(\tau; f, \varphi)$ and $\Phi(z, z'; f_0, \varphi)$ respectively. 
\end{remark}
\begin{proof}
  For the first identity in \eqref{eq:lift-eq}, note both sides are holomorphic by \eqref{eq:RL-commute}, and eigenforms of Hecke operators $T_\pf$ with the same eigenvalue for all but finitely many primes $\pf$.
By strong  multiplicity one, they generate isomorphic automorphic representations. 
Furthermore, the right hand side is a newform of level $N$, whereas the left hand side also has level $N$. So both sides agree up to scalar. 
For the second equality, write the left hand side as linear combinations of (scales of) newforms. Both sides have the same eigenvalue under $T_p$ for almost all primes that split in $F$. Since $(D, N) = 1$, the field $F$ is fixed by the Galois representations associated to the newforms appearing on the left and right hand sides. Applying Chebotarev density theorem and Lemma \ref{lemma:wt1-id} tells us that these newforms are all the same, namely $f_0$.

Finally, Equation \eqref{eq:rr0} is a direct consequence of \eqref{eq:lift-eq} as
  $$
	r_0 \| f_0\|^2 =
  \langle
  I(g; f, \varphi),
  f_0(g)
  \rangle
  =   \langle
f(h),
  \Phi^{1, 1}(h; f_0, \varphi)  
  \rangle
=   \overline{r} \| f\|^2.  \eqno\qedhere
  $$
\end{proof}

\subsection{Unfolding I}
\label{subsec:unfoldI}
We compute the constant $r = r(\varphi)$ from Proposition \eqref{prop:rr0} with $\varphi$ defined in \eqref{eq:varphip}.
The procedure is standard via changing the model of the Weil representation.
A similar case can be found in \textsection 4.2 of \cite{BLY}.

Starting with \eqref{eq:linearize} and \eqref{eq:Poisson}, for $h \in H_1(\Ab)$ and $\varphi^\infty \in \Sc(\hat V)$, we have
\begin{align*}
  \Phi^{1, 1}(h; f_0, \varphi^\infty)
  &=
    \int_{[G]}  \fs_0(g
    ) \sum_{\eta \in \Qb^2, \nu \in V_0}
\Fc(\omega_0(g)\omega(h)\varphi)(\eta \cdot g, \nu)dg ,
\end{align*}
where $\varphi := \varphi^\infty \otimes  \varphi^{1, 1}_\infty \in \Sc(V(\Ab))$.
Since there exists $\gamma \in H_1(\Qb), h' \in H_1(\Ab)$ with $h'_p = 1$ such that $h = \gamma h'$, 
we have
$$
\Fc(\omega_0(g)\omega(h)\varphi)(0, 0, \nu) =
\Fc(\omega_0(g)\omega(h')\varphi)(0, 0, \nu) = 0
$$
for all $g, \nu$ by \eqref{eq:PFT-varphip}.
Therefore the term $\eta = (0, 0)$ in the summation above does not contribute.
The other terms is the $(N(\Qb) \backslash G(\Qb))$-orbit of $(0, 1)$. Writing $g = n(b)m(a) \kappa k_\theta$ and unfolding this gives us
\begin{align*}
  \Phi^{1, 1}(h; f_0, \varphi^\infty)
&    =
\frac12  \sum_{\nu \in V_0}
  \int_{(\Qb\backslash\Ab) \times \Ab^\times}
  \psi(bQ_0(\nu))|a|\chi_F(a)\\
&\times \int_{K}  \fs_0(n(b) m(a) \kappa
  ) 
  \Fc(\omega(\kappa)\omega(h)\varphi)(0, a^{-1}, a \nu)
 d\kappa  db \frac{d^\times a}{|a|^2}.
\end{align*}
The factor $\frac12$ is the volume of $\mathrm{PSO}_2(\Rb)$.


Suppose from now on $h = h_z \in H_1(\Rb)$. 
Then it is not involved in the integral over $\kappa \in K \subset G(\hat\Qb)$. 
Using $\fs_0(g\kappa) = \chi_0(\kappa) \fs_0(g)$ and Lemma \ref{lemma:varphi-cond}, we can write
\begin{align*}
\int_{K} & \fs_0(n(b) m(a) \kappa ) 
  \Fc(\omega(\kappa)\varphi)(0, a^{-1}, a\nu)
  d\kappa\\
  &=
\vol(K_0(N))
    \sum_{\kappa \in K/K_0(N)}
  \fs_0(n(b) m(a) \kappa ) 
    \Fc(\omega(\kappa)\varphi)(0, a^{-1}, a\nu)    \\
  &=
\vol(K_0(N))
    \sum_{N' \mid N}\sum_{0 \le j \le N'-1}
      \fs_0(n(b) m(a) (n(j)w)_{N'} ) 
  \Fc(\omega_0((n(j)w)_{N'})\varphi)((0, a^{-1})w_{N'}, a\nu)    .
\end{align*}
Here we have parametrized $K/K_0(N)$ by\footnote{For $\gamma \in G(\Qb)$ and $M \in \Nb$, we use $\gamma_M \in G(\hat\Qb)$ to denote the element, whose component at $p$ is $\gamma$ if $p \mid M$ and $1$ otherwise.}
$\{(n(j)w)_{N'}: N' \mid N, 0 \le j \le N'-1\}$ as  $N$ is square-free.
If $N' > 1$, i.e.\ there exists $p \mid N'$,
then $((0, a^{-1})w, \nu) = (-a^{-1}, 0, \nu)$ is outside the support of $\Fc_p(\omega_{0, p}(n(j)w)\varphi_p)$ by \eqref{eq:PFT-varphip}. Therefore the corresponding summand vanishes, and only $N'=1$ contributes.
That means
\begin{align*}
  & \frac{\Phi^{1, 1}(h; f_0, \varphi^\infty)}
  {\vol(K_0(N))}
    =
\frac12 \sum_{\nu \in V_0}
\int_{\Ab^\times}  \int_{\Qb \backslash\Ab}
  \psi(bQ_0(\nu))
\fs_0(n(b) m(a)) db
  \Fc(\omega(h)\varphi)(0, a^{-1}, a\nu)
\chi_F(a)
  \frac{d a}{|a|^2}  \\
  & \quad =
  \sum_{\nu \in V_0, \alpha \in \Qb_{>0}}
\int_{ \Rb_{>0}}
\int_{\Qb \backslash\Ab}
    \psi(b Q_0(\nu))
\fs_0(n(\alpha^{-2} b) m(r) ) db
  \Fc(\omega(h)\varphi)(0, (\alpha r )^{-1}, \alpha r \nu)
    \frac{d r }{r^2}  \\
    & \quad =
  \sum_{\substack{\nu \in \df_F^{-1}\\ \alpha^{-1} \in \Nb - DN \Nb\\ \alpha \nu \in \df_F^{-1}}}
    c_{-\alpha^2 Q_0(\nu)}(f_0)
\overline{\chi^\flat_0(-\alpha)}\chi^\flat_F(\alpha)
    \int_{ \Rb_{>0}}
r    e^{2\pi \alpha^2 Q_0(\nu) r^2}
    \Fc_\infty(\omega(h)\varphi_\infty^{1, 1})(0, (\alpha r )^{-1}, \alpha r\nu)
  \frac{d r }{r^2 }    
\end{align*}
In the second step, we have used $\Ab^\times = \Qb^\times (\hat\Zb^\times \Rb_{> 0})$, Remark \ref{rmk:PFT0} and Equation \eqref{eq:f0-adel} to obtain invariance in the integral over $\hat\Zb^\times$.
Since $\alpha$ and $-\alpha$ give the same contribution, we receive a factor of 2. 
In the last step, we applied Equation \eqref{eq:cmf0} and Remark \ref{rmk:PFT0}.

It was computed in \cite[(4.18), (4.19)]{BLY} that
\begin{align*}
  \Fc_\infty(\omega(h)\varphi_\infty^{1, 1})(0, r, \nu)
  &=
\sqrt{y_1y_2} \ebf(r(x_1\nu' + x_2 \nu)) \Fc_\infty(\varphi_\infty^{1, 1})(0, r\sqrt{y_1y_2}, \nu\sqrt{y_2/y_1})\\
  &=
(ry_1y_2 + \nu y_2 + \nu' y_1)
    \ebf(r(x_1\nu' + x_2 \nu))
e^{-\pi (r^2 y_1y_2 + \nu^2 y_2/y_1 + (\nu')^2 y_1/y_2)}. 
\end{align*}
Substituting in this expression gives us the following Fourier expansion for $h \in H_1(\Rb)$
\begin{equation}
  \label{eq:Phi-FE}
  \begin{split}
        \frac{\Phi^{1, 1}(h; f_0, \varphi^\infty)}
  {\vol(K_0(N))}
&=
\sum_{\nu \in \df_F^{-1}, \nu \gg 0}
  \ebf(\nu' z_1 + \nu z_2)
  \sum_{d \mid \sqrt{D}\nu,~ \gcd(d, DN) = 1}
(\chi_0\chi_F)(-d)
  c_{D\Nm(\nu/d)}(f_0) .
  \end{split}
\end{equation}
From this, we obtain the following result.

\begin{proposition}
  \label{prop:r}
  Let
  $r = r(\varphi^\infty)$ be   the constant in Proposition \ref{prop:rr0}.
  Then $r \in \Qb(f_0)^\times$. 
\end{proposition}
\begin{proof}
  Since $D$ and $N$ are co-prime, we can find a prime $\pf$ of $F$ in same narrow class as $\df$ such that $p \nmid DN$ and $c_p(f_0) \neq 0$, where $p = \Nm(\pf)$. Writing $\pf = \df (\nu)$ with $\nu \in F^\times$ totally positive, we obtain
  $$
r\cdot  c_{\nu, \Oc}(  f)
=  c_{\nu, \Oc}(\Phi^{1, 1})
= c_p(f_0) \in \Qb(f_0)^\times.
  $$
Since $f$ is normalized with $c_{1, \Oc}(f) \in \Qb(f)$ and $\Qb(f) \subset \Qb(f_0)$, we have $r \in \Qb(f_0)^\times$. 
\end{proof}

\subsection{Unfolding II}
Now, we will take $f$ to be the base change of $f_0$, and unfold the theta integral $I(g; f, \varphi)$ from \eqref{eq:I}.

For $x \in V(\Qb)$, recall that $H_{1, x} \subset H_1$ is the subgroup stabilizing $x$ defined in Section \ref{subsec:wt-cycle}. If $x = 0$, this is just $H_1$.
If $Q(x) = 0$ but $x \neq 0$, then let $\gamma \in G_F(\Qb)$ be such that $x = \gamma \cdot \smat{a}000$ for some $a \in \Qb^\times$, and 
\begin{equation}
  \label{eq:Hx0}
  B\cap H_1 \cong H_{1, x},~ g \mapsto \gamma g \gamma^{-1}. 
\end{equation}
If $m = -Q(x) \in \Nm(F^\times)$, then there exists $\gamma \in G_F(\Qb) \subset \tH(\Qb)$ such that $x = \gamma \cdot
\smat{}{\sqrt{D}}{-\sqrt{D}}{}$. 
In that case, 
\begin{equation}
  \label{eq:Hx1}
  \SL_2 \cong H_{1, x},~ g \mapsto \gamma g \gamma^{-1}. 
\end{equation}
If $m \not\in \Nm(F^\times)$, then $H_{1, x/F} \cong \SL_{2/F}$ . 

Suppose $x \in V(\Qb)$ is a positive vector.
For any function $\phi: [H_1] \to \Cb$, we define its cycle integral associated with $H_{1,x}$ by
\begin{equation}
  \label{eq:sigmax-f}
  \sigma_{x}(\phi) :=   \int_{[H_{1,x}]}    \phi(h') dh'  
\end{equation}
along the Shimura subvariety obtained from $H_{1,x} \subset H_1$.
It is easy to check that
\begin{equation}
  \label{eq:sigmaxx'}
  \sigma_{\gamma \cdot x}(\phi) =
  \sigma_{ x}(\rho(\gamma^{-1})\phi) 
\end{equation}
for any $\gamma \in H_1(\Qb)$ and $\phi$. 
First, we state a few lemmas needed for the unfolding process.
\begin{lemma}
  \label{lemma:hinf-indep}
  Recall $w_j \in \tH(\Rb)$ from \eqref{eq:wi} and denote $f_j := \rho(w_j) f$.
For any $x \in V(\Qb)$ with $m = Q(x) \neq 0$,   the following functions on $H_1(\Rb)$
  \begin{equation}
    \label{eq:Sigmaf}
    \Sigma_{x, 1}(f, h) := \frac{    \overline{(h\cdot Z^+, x)}}{2\sqrt{D|m|}}
    \sigma_x(\rho(h) f_1),~
    \Sigma_{x, 2}(f, h) :=
    \frac{{(h\cdot Z^+, x)}}{2\sqrt{D|m|}}
    \sigma_x(\rho(h) f_2)
  \end{equation}
  are constant, which we denote by $\Sigma_{x, j}(f)$. 
  If $m < 0$, this is identically 0.
  Otherwise, we have
  \begin{equation}
    \label{eq:Sigmaf-rel}
    \Sigma_{x, j}(f)    =
    \sigma_x(\rho(h_x) f_j)
=\int_{[H_x]} f_j(h'h_x) dh',
\end{equation}
where $h_x \in H_1(\Rb)$ satisfies $x = h_x \cdot (\sqrt{m/D} \smat{}{1}{-1}{})$. 
\end{lemma}
\begin{proof}
  Denote $\hf_{1,x} \subset \hf_1  = \mathfrak{sl}_2^2$ the Lie algebras of $H_{1,x}(\Rb) \subset H_1(\Rb) =\SL_2(\Rb)^2$.
  Choose $h_\infty \in H_1(\Rb)$ such that $h_\infty^{-1} \cdot x =  \sqrt{|m|/D} \smat{0}{1}{-\sgn(m)}0$.
Then $\mathfrak{sl}_2 \cong \hf_{1,x} \subset \mathfrak{sl}_2^2$ via $A \mapsto h_\infty (A, A)h_\infty^{-1}$. 
It is easy to check that $\Sigma_{x, j}(f, h_\infty)$ is right $K_\infty = \SO_2(\Rb)^2$-invariant.
Since $f$ is holomorphic, Equation \eqref{eq:wj-LR} gives us $R_j \Sigma_{x, j} = L_{3-j} \Sigma_{x, j} = 0$  for $j = 1, 2$, where $L_j, R_j$
are the raising and lowering operators in $\slf_{2, \Cb}^2= \hf_{1, \Cb}$ defined in \eqref{eq:LR}.
Since $\Sigma_{x, j}$ is left $H_{1,x}(\Rb)$-invariant, it is annihilated by $\hf_{1, x, \Cb}, R_j, L_{3-j}, \mathfrak{so}_{2, \Cb}^2$, which generate $\hf_{1, \Cb}$. Therefore it is independent of $H_1(\Rb)$.

For the second claim, notice
$$
(h \cdot Z^+, x) = (Z^+, h^{-1} \cdot x) =  (1 + \sgn(m))  \sqrt{D|m|}, 
$$
which is zero for $m < 0$. 
The last claim follows easily after substitution. 
\end{proof}

\begin{lemma}
  \label{lemma:arch-int}
  For any $x \in V(\Rb)$ with $Q(x) = m \neq 0$, we have
  \begin{equation}
    \label{eq:arch-int}
    \int_{H_x(\Rb) \backslash H_1(\Rb)} e^{-2\pi vQ_+(h^{-1}x)} dh
    = \frac1{mv} e^{-2\pi mv}.
  \end{equation}
\end{lemma}

\begin{proof}
  After translation by $h$, we can suppose $x = \sqrt{|m|/D}\smat{}{1}{-\sgn(m)}{}$. Then  $\SL_2(\Rb) \cong H_x(\Rb) \backslash H_1(\Rb)$ via $g \mapsto (g, 1)$. The measure $dh = dg$ is $d\mu(z) \frac{d\theta}{2\pi}$ when $g = g_z\kappa_\theta$.
  Direct calculation gives
  $$
2\pi v Q_+(g_z\kappa_\theta x) = \frac{\pi v m}y (|z|^2 + 1). 
$$
Therefore, we have 
\begin{align*}
  \int_{H_x(\Rb) \backslash H_1(\Rb)} e^{-2\pi vQ_+(h^{-1}x)} dh
  &=
    \int_\Hb e^{-\frac{\pi v m}y (|z|^2 + 1)} d\mu(z),
\end{align*}
which yields \eqref{eq:arch-int} after standard calculations. 
\end{proof}
Now we are ready to start the unfolding process. 
Take $g_\infty = g_\tau = n(u)m(\sqrt{v})\in \SL_2(\Rb), g_0 = g_\infty \smat{1}{}{}{-1} \in G(\Rb) \subset G(\Ab)$, and $h_0 = w_j \in H(\Rb)$. Then $\nu(h_0) = \nu(g_0)$ and
\begin{align*}
&v^{-1}  I(g_0; f, \varphi)
  =
  v^{-1}
  \int_{[H_1]}  f( hw_j)
    \overline{\theta^{1, 1}(g_\infty \smat{1}{}{}{-1},  hw_j; \varphi)} dh\\
  & = \sum_{x \in H_1(\Qb)\backslash V(\Qb)}
\ebf(-Q(x)u)
    \int_{H_{1,x}(\Ab)\backslash H_1(\Ab)} 
    \overline{    \varphi((h^\infty)^{-1}x)(\cL(w_j)\varphi_\infty^{1, 1})(\sqrt{v} h_\infty^{-1}\cdot x)}
    \sigma_x(\rho(h) f_j) dh^\infty dh_\infty.
\end{align*}
When $x = 0$, the function $\varphi_\infty^{1, 1}$, hence integrand, vanishes identically. On $V(\Qb) - \{0\}$, the group $H_1(\Qb) = \mathrm{Spin}_V(\Qb)$ acts transitively.
For each $m \in \Qb$, we then choose $x_m \in V(\Qb)$ with $Q(x_m) = m$.
As in Section \ref{subsec:wt-cycle}, we can find a finite subset $S(\varphi, x_m) \subset   H(\hat\Qb)$ satisfying \eqref{eq:Svarphix0} and
 rewrite the unfolding as
\begin{align*}
v^{-1}    I(g_0; f, \varphi)
  & =
\vol(K)
    \sum_{\substack{m \in \Qb\\s \in S(\varphi, x_m)}}
\ebf(-mu)\overline{    \varphi(s^{-1}\cdot x)}\\
& \quad\times    \int_{H_{x_m}(\Rb)\backslash H_1(\Rb)} 
    \overline{ (\cL(w_j)\varphi_\infty^{1, 1})(\sqrt{v} h_\infty^{-1}x_m)}
    \sigma_{x_m}(\rho(s, h_\infty) f_j) dh_\infty.
\end{align*}
To simplify the integral above, we can apply \eqref{eq:varphi-w} and Lemmas \ref{lemma:hinf-indep}, \ref{lemma:arch-int} to obtain
\begin{align*}
      \int_{H_{x_m}(\Rb)\backslash H_1(\Rb)} 
&    \overline{ (\cL(w_j)\varphi_\infty^{1, 1})(\sqrt{v} h_\infty^{-1}x_m)}
    \sigma_{x_m}(\rho(s, h_\infty) f_1) dh_\infty
\\
&=      \int_{H_{x_m}(\Rb)\backslash H_1(\Rb)} 
    \overline{ \varphi_\infty^{-1, -(-1)^j}(\sqrt{v} x_m)}
  \sigma_{x_m}(\rho(s, h_\infty) f_1) dh\\
&=
      \int_{H_{x_m}(\Rb)\backslash H_1(\Rb)} 
  i {D^{-1}(h_\infty^{-1}\cdot x_m, X^+ + (-1)^j i Y^+)}e^{-4\pi vQ_+(h_\infty^{-1}x_m)}
  \sigma_x(\rho(s, h_\infty) f) dh\\  
&=
(-1)^{j+1} \sqrt{D|m|} \Sigma_{x_m, j}(\rho(s)f)
  \int_{H_{x_m}(\Rb)\backslash H_1(\Rb)} 
  e^{-2\pi vQ_+(h_\infty^{-1}x_m)} dh\\
&  = \frac{(-1)^{j+1}\sqrt{D}}{\sqrt{m} v} \Sigma_{x_m, j}(\rho(s)f) e^{-2\pi m v}. 
\end{align*}
By \eqref{eq:sigmaxx'} and \eqref{eq:Svx}, the following sum
\begin{equation}
  \label{eq:Cmv}
  \begin{split}
      \Cc_{m, j}(f; \varphi)
&  :=
  \sqrt{D/|m|}
  \sum_{s \in S(\varphi, x_m)}
\overline{\varphi(s^{-1} x_m)}
  \Sigma_{x_m, j}(\rho(s)f)
  \end{split}
\end{equation}
is independent of the choice of $x_m$.
Putting these together gives us the following result.
\begin{theorem}
  \label{thm:unfold2}
  In the notations above, we have
  $$
 r_0(\varphi) f_0(-\overline{\tau})
=  I(g_0; f, \varphi) =
\sum_{m \in \Qb_{> 0}} (-1)^{j+1}
 \Cc_{m, j}(f; \varphi)
 \ebf(-m \overline\tau)
  $$
  for $j = 1, 2$.
  In particular, we have $r_0(\varphi) = \Cc_{1, 1}(f; \varphi) = -\Cc_{1, 2}(f; \varphi)$. 
\end{theorem}

\begin{proof}
  The first equality is simply Proposition \ref{prop:rr0}, plus the relation  $f_0^\#(g_0) = f_0(-\overline{\tau})$. 
  The second equality follows from the deductions preceding Theorem \ref{thm:unfold2}. 
\end{proof}


\section{Motivic action conjecture}

In this section, we deduce from Theorem~\ref{thm:unfold2} the motivic action conjecture for Doi--Nagunuma lifts of weight one modular forms  (Theorem~\ref{thm:motivic_action}).

\subsection{Statement of the conjecture}

We start by briefly recalling the statement of the motivic action conjecture for Hilbert modular forms of weight one over real quadratic fields $F/\Q$, formulated in \cite{Horawa}. We keep the notation of Section~\ref{sec:prelim}.

Let $f$ be an adelic Hilbert modular form for $G_F$ of weight $(1,1)$ and level $K_1(N)$. Assume that $f$ is a normalized eigenform and write $E := \Q(f) \subseteq \C$ for the number field generated by the Hecke eigenvalues of $f$. Associated to a level $K = K_1(N)$, there is a Shimura variety $X := X_K$ defined over $\Q$, together with a Hodge bundle $\omega$, such that $f$ defines an $E$-rational section $[f] \in H^0(X, \omega) \otimes E$. In fact, $H^0(X, \omega)$ has an action of the Hecke operators, and if we let $H^0(X, \omega)_f$ be the $f$-isotypic component of $H^0(X, \omega) \otimes_\Q E$, then $H^0(X, \omega)_f$ is 1-dimensional, spanned by $[f]$.  

The point of the motivic action conjecture is to consider the higher coherent cohomology group $H^1(X, \omega)$. Since the variety $X$ is not proper, we consider a suitable toroidal compactification $X^\tor$ defined over $\Q$. The line bundle $\omega$ admits two natural extensions $\omega^{\rm{sub}}$ and  $\omega^{\rm{can}}$ to $X^\tor$, and we define the $\Q$-vector space:
\begin{equation}
	H^1(X, \omega) := \mathrm{Im}( H^1(X, \omega^{\rm{sub}}) \to  H^1(X, \omega^{\rm{can}})),
\end{equation}
as in~\cite{Harris_Periodinvariants, Horawa}. It again admits a Hecke action and we consider the $f$-isotypic component $H^1(X, \omega)_f$, an $E$-vector space. Associated with $w_j \in G_F(\R)$ (defined in \eqref{eq:wi}), there are automorphic forms $f_j(g) := f(g \cdot w_j)$ for $j = 1,2$, which give rise to a natural basis of this cohomology group over $\C$:
\begin{equation}\label{eqn:omega_f^j}
	\omega_f^j := [f_j] \in H^1(X_\C, \omega)_f \iso H^1(X, \omega)_f \otimes_E \C
\end{equation}
(c.f.\ \cite{Harris_Periodinvariants}, where $\omega$ is denoted $\mathcal E_{\underline 1, 1}$). Very roughly, the motivic action conjecture now predicts that the rationality of these classes is related to certain Stark units.

Let $\varrho_f \colon \Gal_{F} \to \GL_2(\C)$ be the Artin representation associated with $f$~\cite{Deligne_Serre, Rogawski_Tunnell}. This representation is odd and clearly has traces in $E$, and this means it can be defined over $E$. Then the trace zero adjoint representation $\Ad^0 \varrho \colon \Gal_F \to \GL_3(E)$ is obtained from the conjugation action of $\Gal_F$ on traceless matrices. Since $\Ad^0 \varrho_f$ factors through a finite quotient, there is a finite Galois extension $L/F$ such that
\begin{equation}
	\Ad^0 \varrho \colon \Gal(L/F) \to \GL_3(E).
\end{equation}
(In fact, $\Ad^0 \varrho$ might be defined over a smaller field than the field $E$ of Fourier coefficients of $f$, but we extend it here to $E$ for convenience.)

\begin{definition}
	The {\em Stark unit group} associated with $f$ is
	$$U_f := \Hom_{E[\Gal(L/F)]}(\Ad^0 \varrho, \O_L^\times \otimes_\Z E).$$
\end{definition}

Given the two embeddings $F \hookrightarrow \R$, we extend them to embeddings $\tau_1, \tau_2 \colon L \hookrightarrow \C$, and get two associated complex conjugation automorphisms $c_1, c_2 \in \Gal(L/F)$. 

\begin{proposition}[{\cite[Cor.\ 2.13]{Horawa}}]
	There is a natural decomposition
	$$U_f \iso ((\Ad^0 \varrho_f)^{c_1})^\vee \oplus ((\Ad^0 \varrho_f)^{c_2})^\vee$$
	into 1-dimensional vector spaces $((\Ad^0 \varrho_f)^{c_i})^\vee$.
\end{proposition}

For $i = 1, 2$, consider a vector $v_i \in \Ad^0 \varrho_f$ such that $\varrho_f(c_i) v_i = v_i$, well-defined up to $E^\times$. Moreover, for each $\varphi_j \in U_f$ lying in the line $((\Ad^0 \varrho_f)^{c_1})^\vee$ in the above decomposition, we get a unit:
$$u_{ij} = \varphi_j(v_i) \in \O_{L}^\times \otimes_\Z E.$$

We obtain the matrix of logarithms of these units $R_f := (\log|u_{ij}|)_{i,j}$ and Stark's conjecture predicts that
$$L(\Ad^0 \varrho_f, 1) \sim_{E^\times} \frac{\pi^{2d}}{f_\varrho^{1/2}} \det R_f,$$
where $f_\varrho$ is the conductor of $\Ad^0 \varrho_f$ (c.f.\ \cite[Prop.\ 2.17]{Horawa}).

The motivic action conjecture is then the following.

\begin{conjecture}[{ \cite[Conjecture 4.17]{Horawa} }]\label{conj:motivic_action}
	An $E$-basis of $H^1(X, \omega)_f$ is given by 
	$$R_f^{-1} \cdot \begin{pmatrix}
		\omega_f^1 \\
		\omega_f^2
	\end{pmatrix}.$$
Explicitly, the basis is:
	\begin{align}
		\frac{\log|u_{22}| \omega_f^1 - \log|u_{12}| \omega_f^2}{\det R_f}, \\
		\frac{-\log|u_{21}| \omega_f^1  + \log|u_{11}| \omega_f^2}{\det R_f}.
	\end{align}
\end{conjecture}

\subsection{Conjecture for base change forms}

Now, suppose that $f_0$ is a weight one modular form and $f$ is a base change of $f_0$, i.e.\ $\varrho_f = \varrho_{f_0}|_{G_F}$. Then the unit group $U_f$ is
$$U_f \iso \Hom(\Ad^0 \varrho_0, \O_L^\times \otimes E) \oplus \Hom(\Ad^0 \varrho \otimes \chi_F, \O_L^\times \otimes E).$$
We fix $\varphi_0 \in \Hom(\Ad^0 \varrho_0, \O_L^\times \otimes E)$ and $\varphi_0^F \in \Hom(\Ad^0 \varrho_0 \otimes \chi_F, \O_L^\times \otimes E)$, together with $v_0 \in (\Ad^0 \varrho_0)^c$, and define
$$u_{f_0} := \varphi_0(v_0),\ u_{f_0}^F := \varphi_0^F(v_0).$$
One can then check that:
\begin{equation}
	R_f = \begin{pmatrix}
		1 & -1 \\
		1 & 1
	\end{pmatrix}  \begin{pmatrix}
		\log|u_{f_0}| & 0 \\
		0 & \log|u_{f_0}^F|
	\end{pmatrix} \begin{pmatrix}
		1 & -1 \\
		1 & 1
	\end{pmatrix}^{-1}
\end{equation}
(c.f.\ \cite[Corollary 6.5]{Horawa}). In particular, Conjecture~\ref{conj:motivic_action} has the following form.

\begin{conjecture}[{\cite[Conjecture 6.7]{Horawa}}]\label{conj:motivic_action_BC}
	An $E$-basis of $H^1(X_\Q, \omega)_{f}$ is given by:
	$$  \frac{\omega_f^1 + \omega_f^2}{\log|u_{f_0}|}, \frac{\omega_f^1 - \omega_f^2}{\log|u_{f_0}^F|}.$$
\end{conjecture} 
Our goal is to prove this conjecture. 

\begin{theorem}\label{thm:motivic_action}
	Suppose $f_0$ is a modular form of square-free level $N$ co-prime to $2 D$, and let $f$ be the base change of $f_0$ to $F = \Q(\sqrt{D})$. Then:
	\begin{enumerate}
		\item An $E$-basis of $H^1(X_\Q, \omega)_{f}$ is given by:
		$$  \frac{\omega_f^1 + \omega_f^2}{\langle f_0, f_0 \rangle}, \frac{\omega_f^1 - \omega_f^2}{ \frac{\langle f, f \rangle}{\langle f_0, f_0 \rangle}}.$$
		\item If the projective image of $\varrho_0$ is $A_5$, assume Stark's Conjecture~\cite{Stark-II} for $\Ad^0 \varrho$ and $\Ad^0 \varrho \otimes \chi_F$. Then Conjecture~\ref{conj:motivic_action_BC} is true.
	\end{enumerate}
\end{theorem}

\subsection{Proof of Theorem~\ref{thm:motivic_action}}

The proof of the theorem will proceed in 3 steps:
\begin{enumerate}
	\item There is a rational involution $s \colon X \to X$ (the swap map) which allows to decompose the 2-dimensional $E$-vector space $H^1(X, \omega)_f$ into two 1-dimensional $E$-vector spaces $H^1(X, \omega)_f^\pm$ such that
	\begin{align}
		\eta^+ := \frac{\omega_f^1 + \omega_f^2}{\langle f_0, f_0 \rangle} \in H^1(X, \omega)_f^+ \otimes \C,  \\
		\eta^- := \frac{\omega_f^1 - \omega_f^2}{\frac{\langle f, f \rangle}{\langle f_0, f_0 \rangle}} \in H^1(X, \omega)_f^- \otimes \C.
	\end{align}
	The goal is then to prove that $\eta^+$ and $\eta^-$ are both rational.
	\item Theorem~\ref{thm:unfold2} implies that $\eta^-$ is rational.
	\item One then uses Serre duality (Proposition~\ref{prop:Serre_duality_consequence}) to check that $\eta^+$ is rational.
\end{enumerate}

We start by introducing the swap map. For any $\Q$-algebra $R$, we have a map $G_F(R) \to G_F(R)$ induced by the non-trivial Galois element $\sigma \in \Gal(F/\Q)$ via the following diagram:
\begin{center}
	\begin{tikzcd}[column sep = small]
		\Spec(R \otimes_\Q F) \ar[rd,swap,  "g"] \ar[rr, "(1 \otimes \sigma)^\ast"] & & \Spec(R \otimes_\Q F) \ar[ld, "s(g)"] \\
		&  \GL_{2, F}
	\end{tikzcd}
\end{center}
This gives a map $s \colon G_F \to G_F$ of algebraic groups, which induces a map at the level of Shimura varieties.

\begin{proposition}\label{prop:swap}
	Let $X$ be the open Hilbert modular surface of level $\Gamma_1(N \O_F)$ for an integer $N$. Then the swap map $s$ induces a map $s \colon X \to X$ which defined over $\Q$. Moreover, there is a choice of toroidal compactification $X \hookrightarrow X^\tor$ such that $s$ extends to a map $s \colon X^\tor \to X^\tor$ defined over $\Q$. Moreover, $s^\ast \omega = \omega$ and hence $s^\ast \omega^{\mathrm{sub}} \iso \omega^{\mathrm{sub}}$ and $s^\ast \omega^{\mathrm{can}} \iso \omega^{\mathrm{can}}$.
\end{proposition}
\begin{proof}
	First, note that $s(\Gamma_1(N \O_F)) = \Gamma_1(N \O_F)$, so $s$ induces a map $X_\C \to X_\C$. As in~\cite[\S 5]{Deligne_Ribet}, recall that $X$ is a moduli space of polarized abelian varieties with $\O_F$-multiplication with $\Gamma_0(N \O_F)$-level structure, i.e.\ an $\O$-linear immersion $\alpha_N \colon \O \otimes \mu_N \hookrightarrow A$. Then the map $s$ induces the following map on moduli:
	$$(A/S, \iota, \lambda, \alpha_N) \mapsto (A/S, \iota^\sigma, \lambda, \alpha_N^\sigma),$$
	where $\iota^\sigma(x) = \iota(x^\sigma)$ and $\alpha_N^\sigma(x \otimes \zeta) = \alpha_N(x^\sigma \otimes \zeta)$. Since this is clearly rational, the map $s$ induces a rational map $s \colon X \to X$. 
	
	Consider the minimal compactification of $X(\C) \hookrightarrow X(\C)^\ast$, obtained by adding the finitely many cusps $\Gamma_1(N\O_F) \backslash \P^1(F)$. The action of $s$ extends to the natural Galois action of $\sigma$ on $\Gamma_1(N\O_F) \backslash \P^1(F)$. For each $c \in \Gamma_1(N\O_F) \backslash \P^1(F)$:
	\begin{enumerate}
		\item if $c^\sigma = c$, then choose any admissible polyhedral cone decomposition for $c$,
		\item if $c^\sigma \neq c$, then choose any admissible polyhedral cone decomposition $\Sigma$ for $c$, and then the associated polyhedral cone decomposition $\Sigma^\sigma$ for $c^\sigma$.
		.	\end{enumerate}
	This gives rise to a toroidal compactification
	\begin{center}
		\begin{tikzcd}
			& X^\tor \ar[d, "\pi"] \\
			X \ar[r, hook] \ar[ru, hook] & X^\ast
		\end{tikzcd}
	\end{center}
	We extend $s \colon X \to X$ to a map $s \colon X^\tor \to X^\tor$ by defining:
	\begin{enumerate}
		\item if $c^\sigma = c$, then $s|_{\pi^{-1}(c)} = 1_{\pi^{-1}(c)}$,
		\item if $c^\sigma \neq c$, then $s|_{\pi^{-1}(c)} \colon \pi^{-1}(c) \to \pi^{-1}(c^\sigma)$ is the natural map induced by $\sigma$.
	\end{enumerate}
	This gives a rational map $s \colon X^\tor \to X^\tor$, as required. In that notation of~\cite{Harris_Periodinvariants}, it is easy to check that $s^\ast \mathcal E_{(k_1, k_2; r)} = \mathcal E_{(k_2, k_1; r)}$, and $s^\ast \omega = \omega$ follows from the identification $\omega \iso \mathcal E_{(1,1; 1)}$.
\end{proof}

We next observe that the swap map switches $\omega_f^1$ and $\omega_f^2$ for a base change form $f$. 

\begin{lemma}
	In general, $s^\ast \omega_f^1 = \omega_{s^\ast f}^2$. Moreover, for a base change form $f$, $s^\ast f = f$, so $s^\ast \omega_f^1 = \omega_f^2$.
\end{lemma}
\begin{proof}
	Note that $s_\R(w_1) = w_2$ under the map $s_\R \colon G(\R) \to G(\R)$. Therefore, $s^\ast [f_1] = [s^\ast f_1]$ and $s^\ast f_1(g) = f(s(w_1 g)) = f(w_2 s(g)) = (s^\ast f)_2(g)$. Finally, if $f$ is a base change form, then by Proposition~\ref{prop:rr0}, $\Phi^{1,1}(h, f_0, \varphi) = r \cdot f(h)$ and we note that $s^\ast \Phi^{1,1}(h, f_0, \varphi) = \Phi^{1,1}(h, f_0, \varphi)$ by \eqref{eq:Phi}.
\end{proof}

Since $s^\ast$ defines an involution on $H^1(X, \omega)_f$ for a base change form $f$, we can define:
\begin{equation}
	H^1(X, \omega)_f^\pm = \{\eta \in H^1(X, \omega)_f \ | \ s^\ast \eta = \pm \eta \}.
\end{equation}

\begin{corollary}
	The space $H^1(X, \omega)_f^\pm \otimes \C$ is spanned by $\omega_f^1 \pm \omega_f^2$. In particular, $$\dim(H^1(X, \omega)_f^\pm) = 1.$$
\end{corollary}

Next, we want to use use Theorem~\ref{thm:unfold2} to get an $E$-rational functional on $H^1(X, \omega)_f^-$. Recall the weighted cycle~\eqref{eq:Z1varphi}:
\begin{align}
	 Z_s & := Z_1(x_0, s, K_1), & s \in S(\varphi, x_0)\\
	 Z_1(\varphi) & := \sqrt{D}  \cdot Z_1(1, \varphi, K_1) := \sqrt{D} \sum_{s \in S(\varphi, x_0)} \varphi(s^{-1} \cdot x_0) Z_s.
\end{align}
whose components $Z_s$ are defined over $E_N := \Q(\zeta_N)$ and whose coefficients are in the field of values of $\sqrt{D} \cdot \varphi$. 

We now fix the choice of $\varphi$ in~\eqref{eq:varphip} (or indeed its complex conjugate which is what we will eventually use). Note that:
\begin{equation}
	\varphi = \prod_{p|N} \frac{1}{\mathfrak g_{\chi_{0, p}}} \cdot \prod_{p|D} \frac{1}{\mathfrak g_{\chi_{F, p}}} \cdot \varphi_0,
\end{equation}
where $\varphi_0$ has values in $E_0 = \Q(\chi_0)$. By definition, $\mathfrak g_{\chi_{0, p}} \in E_{0, N} := E_0E_N$ for all $p|N$. 


The Gauss sums for $p | D$ were computed by Gauss:
\begin{equation}
	\mathfrak g_{\chi_{F, p}} = \begin{cases}
		\sqrt{p} & p \equiv 1 \mod 4, \\
		i\sqrt{p} & p \equiv 3 \mod 4,  \\
	\end{cases}
\end{equation}
Under our assumption that $D$ is odd, there is an even number of primes $p \equiv 3 \mod 4$ dividing $D$, and hence $\sqrt{D} \cdot \prod\limits_{p|D} \frac{1}{\mathfrak g_{F,p}} \in \Q^\times$. Altogether, this shows that $\sqrt{D} \cdot \varphi$ has values in $E_{0,N}$. In summary,
\begin{equation}\label{eqn:Z_1(varphi)_fd_of_def}
	Z_1(\varphi) \in \mathrm{Div}(X_{\Q(\zeta_N)}) \otimes_{\Z} E_{0, N}.
\end{equation}

\begin{lemma}\label{lemma:rat_of_weighted_cycle}
	For $\varphi$ as in~\eqref{eq:varphip}, the divisor $Z_1(\varphi)$ is invariant under the diagonal action of $\Gal(E_{0,N} / E_0) \iso \Gal( E_N/ E_N \cap E_0) \subseteq \Gal(E_N/\Q) \iso (\Z/N\Z)^\times$.
\end{lemma}
\begin{proof}
	Recall that action of $\sigma_\delta \in \Gal(E_N/\Q)$ on the connected components of the cycle is:
	\begin{align}\label{eqn:Galois_action_on_HZ_components}
		Z_1(X, h, K_1)^{\sigma_\delta}  = Z_1(x_0, h \delta^{-1}, K_1)
	\end{align}
	by Equations~\eqref{eq:E-act} and~\eqref{eq:Z1}. 
	
	It is enough to compute the Galois action on $\varphi_p$ for $p | N$. Taking $\delta_a = \begin{pmatrix}
		a & 0 \\
		0 & 1
	\end{pmatrix}$ for $a \in \Z$ representing a class in $(\Z/N\Z)^\times$, we see from~\eqref{eq:varphip} that
	\begin{equation}\label{eqn:delta_a_on_phi_p}
		\varpi_p(\delta_a^{-1} x) = \chi_0(a) \varphi_p(x).
	\end{equation}
	Therefore, for $\sigma_a \in \Gal(E_K/E_0) \subseteq (\Z/N\Z)^\times$, we have
	\begin{align*}
		\varphi_p(\delta_a^{-1} x)^{\sigma_a} & = (\chi_0(a) \varphi_p(x))^{\sigma_a} \\
		& = \chi_0(a) \varphi_p(x)^\sigma.
	\end{align*}
	Since $\mathfrak g_{\chi_0, p}^\sigma = \overline{\chi_0(a)} \mathfrak g_{\chi_0, p}$ and the values of $\varphi_p$ are otherwise $E_0$-rational, this gives:
	\begin{equation}\label{eqn:diag_action_on_varphi_p}
		\varphi_p(\delta_a^{-1} x)^{\sigma_a} = \chi_0(a) \overline{\chi_0(a)} \varphi_p(x) = \varphi_p(x).
	\end{equation}
	Altogether, this gives that:
	\begin{align*}
		Z_1(\varphi)^{\sigma_a} & = \sum_{s \in S(\varphi, x_0)} (\sqrt{D} \varphi(s^{-1} x_0))^{\sigma_a} Z_1(X, s, K_1)^{\sigma_a} \\
		& =  \sum_{s \in S(\varphi, x_0)} \sqrt{D} \varphi(s^{-1} x_0)^{\sigma_a} Z_1(X, s \delta_a^{-1}, K_1) & \text{\eqref{eqn:Galois_action_on_HZ_components}} \\
		& =  \sum_{s' \in S(\varphi, x_0)} \sqrt{D} \varphi(\delta_a^{-1} s'^{-1} x_0)^{\sigma_a} Z_1(X, s', K_1) & s' = s \delta_a^{-1} \\
		& =  \sum_{s' \in S(\varphi, x_0)} \sqrt{D} \varphi(s'^{-1} x_0) Z_1(X, s', K_1) & \eqref{eqn:diag_action_on_varphi_p} \\
		& = Z_1(\varphi),
		\end{align*}
	completing the proof.
\end{proof}

For any $s \in S(\varphi, x_0)$, we write $\iota_s \colon Z_s \hookrightarrow X_{E_N}$ for the embedding, and observe that $\iota^\ast(\omega) \iso \Omega^1_{Z_s}$. Lemma~\ref{lemma:rat_of_weighted_cycle} (applied to $\overline \varphi$ instead of $\varphi$) immediately gives the following corollary.

\begin{corollary}
\label{cor:Cc}
	The functional
	\begin{align*}
		\mathcal C_{\overline \varphi} \colon  H^1(X_\C, \omega)_f & \to \C \\
		\eta & \mapsto  \sqrt{D} \sum_{s \in S(\varphi, x_1)} \overline{\varphi(s^{-1} \cdot x_0)} \int\limits_{Z_1(x_0, s, K_1)} \iota_s^\ast(\eta),
	\end{align*}
	is $E$-rational, i.e.\ it is the base change of a functional
	$$\mathcal C_{\overline \varphi} \colon H^1(X, \omega)_f \to E.$$
\end{corollary}

Finally, Theorem~\ref{thm:unfold2} gives the following rationality statement.

\begin{corollary}\label{cor:eta^-_rational}
	We have that:
	\begin{align}
		\mathcal C_{\overline \varphi}\left( \eta^+ \right) & = 0, \label{eq:Con+=0} \\
		\mathcal C_{\overline \varphi}\left( \eta^-  \right) & \in E^\times.\label{eq:Con-_rational}
	\end{align}
	Therefore: 
	\begin{equation}\label{eq:eta^-_rational}
		\eta^- \in H^1(X, \omega)_f^- \subseteq H^1(X, \omega)_f.
	\end{equation}
\end{corollary}
\begin{proof}
	Note that $\mathcal C_{\overline{\varphi}}(\omega_f^j) = \mathcal C_{1, j}(f; \varphi)$ by definition of $\mathcal C_{1, j}(f; \varphi)$ in~\eqref{eq:Cmv} and Lemma~\ref{lemma:hinf-indep}. Therefore, the two Equations~\eqref{eq:Con+=0} and \eqref{eq:Con-_rational} follow from Theorem~\ref{thm:unfold2} together with Proposition~\ref{prop:rr0} and Proposition~\ref{prop:r}.
	
	Finally, since $H^1(X, \omega)_f^-$ is 1-dimensional and $\eta^- \in H^1(X, \omega)_f^- \otimes \C$, there is a constant $c^- \in \C^\times$ such that $c^- \cdot \eta^-  \in H^1(X, \omega)_f^-$. Then $c^- \cdot \mathcal C_{\overline \varphi}(\eta^-) = \mathcal C_{\overline \varphi}(c^- \cdot \eta^-) \in E^\times$ which together with~\eqref{eq:Con-_rational} gives $c^- \in E^\times$, so $\eta^- \in H^1(X, \omega)_f^-$.
\end{proof}

To complete the proof of Theorem~\ref{thm:motivic_action}, it remains to prove the rationality of $\eta^+$. There is a Serre duality pairing $\langle -, - \rangle_{\rm SD} \colon H^1(X, \omega) \times H^1(X, \omega) \to \Q$ which we can use to prove the following. 

\begin{proposition}[{\cite[Prop.\ 5.14]{Horawa}}]\label{prop:Serre_duality_consequence}
	There is an alternating pairing
	\begin{align*}
		\langle -, - \rangle \colon H^1(X, \omega)_f \times H^1(X, \omega)_f \to \Q(f) 
	\end{align*}
	such that $\langle \omega_f^1, \omega_f^2 \rangle = \langle f, f \rangle$. In particular, the decomposition $H^1(X, \omega)_f = H^1(X, \omega)_f^+ \oplus H^1(X, \omega)_f^-$ is a polarization with respect to this pairing. 
\end{proposition}

\begin{corollary}\label{cor:eta^+_rational}
	We have that
	$$\eta^+ \in H^1(X, \omega)_f^+ \subseteq H^1(X, \omega)_f.$$
\end{corollary}
\begin{proof}
	Recall that $\eta^- \in H^1(X, \omega)_f^-$ by \eqref{eq:eta^-_rational}. Since $H^1(X, \omega)_f^+ \otimes \C$ is 1-dimensional, there is a constant $c^+ \in \C^\times$ such that $c^+ \cdot \eta^+ \in H^1(X, \omega)_f^+$. Therefore, $c^+ \langle \eta^+, \eta^- \rangle \in E^\times$. On the other hand:
	\begin{align*}
		 \langle \eta^+, \eta^- \rangle & = \left\langle \frac{\omega_f^1 + \omega_f^2}{\langle f_0, f_0 \rangle}, \frac{\omega_f^1 - \omega_f^2}{\frac{\langle f, f \rangle}{\langle f_0, f_0 \rangle}} \right\rangle \\
		& = -2
	\end{align*}
	showing that $c^+ \in E^\times$, and hence $\eta^+ \in H^1(X, \omega)_f^+$.
\end{proof}

With these ingredients in place, we are ready to prove Theorem~\ref{thm:motivic_action}.

\begin{proof}[Proof of Theorem~\ref{thm:motivic_action}]
	Part (1) is Equation~\eqref{eq:eta^-_rational} in Corollary~\ref{cor:eta^-_rational} together with Corollary~\ref{cor:eta^+_rational}.
%
	
	For part (2), we first use the well-known relationship between the Petersson norm of a Hilbert modular form $f$ and the adjoint $L$-function. For example, by~\cite[Theorem 7.1]{Hida_Tilouine}\footnote{The factor of $\sqrt{D}$ comes from a difference in normalization of measures. Indeed, \cite{Hida_Tilouine} use the normalizations of~\cite[Section 4]{Hida-p-adic_tot_real}, where the additive measure gives $\O_v$ volume 1 instead of $|\mathfrak d|_v^{1/2}$. C.f.\ also~\cite[Prop.\ 6.6]{Ichino_Prasanna} where the same normalization is used.}:
	\begin{equation}\label{eqn:<f,f>}
		\langle f, f \rangle \sim_{\Q^\times} \pi ^{-2d} \cdot \sqrt{D} \cdot L(f, \Ad,1). 
	\end{equation}
	In fact, there is even a completely explicit formula for the rational constant, but we omit it here.
	This gives:
	\begin{align}
		\langle f, f \rangle & \sim_{\Q^\times} \pi^{-4} \cdot \sqrt{D} \cdot L(f, \Ad, 1) \\
		\langle f_0, f_0 \rangle & \sim_{\Q^\times} \pi^{-2} \cdot L(f_0, \Ad, 1) \label{eqn:<f_0,f_0>}
	\end{align}
	Next, we want to use Stark's conjecture for $\Ad^0 \varrho_{0}$ and $\Ad^0 \varrho_0 \otimes \chi_F$ to connect this to logarithms of units. We need to consider the different possibilities for the projective image of $\varrho_0$:
	\begin{itemize}
		\item Suppose $\varrho_0$ has dihedral projective image, i.e.\ $f_0$ is associated with a character $\chi$ of a quadratic extension $K/\Q$. Then $\Ad^0 \varrho_0 \iso \Ind_{K}^\Q (\chi/\chi^c) \oplus \chi_K$. If $K$ is imaginary quadratic, then the unit group $U_{f_0}$ is associated with the character $\chi/\chi^c$ of $K$ and the unit group $U_{f_0}^F$ is associated with a character of the third intermediate extension between $\Q$ and $KF$; in either  case, the conjecture was proved by Stark~\cite{Stark-IV}. If~$K$ is real quadratic, the conjecture for both $\Ad^0 \varrho_0$ and $\Ad^0 \varrho_0 \otimes \chi_K$ follows from the class number formula.
		
		\item If $\varrho_0$ has solvable projective image which is not dihedral (i.e.\ $A_4$ or $S_4$), then $\Ad^0 \varrho_0$ and $\Ad^0 \varrho_0 \otimes \chi_F$ both have rational traces, and the conjecture was proved by Stark~\cite{Stark-II}.
		
		\item If $\varrho_0$ has projective image $A_5$ (the only non-solvable case), then $\Ad^0 \varrho_0$ does not have rational traces and we have to assume the Stark's conjecture in this case.
	\end{itemize}

	In all cases, under the assumptions of (2), we obtain:
	\begin{align}
		L(1, \Ad^0 \varrho_0)  & \sim_{E^\times} \frac{\pi^2}{f_{\Ad^0 \varrho_0}^{1/2}}  \cdot \log|u_{f_0}|, \label{eqn:Ad_varrho0} \\ 
		L(1, \Ad^0 \varrho_0 \otimes \chi_F)  & \sim_{E^\times} \frac{\pi^2}{f_{\Ad^0 \varrho_0}^{1/2} |D|^{1/2}}  \cdot \log|u_{f_0}^F|. \label{eqn:Ad_varrho0_F}
	\end{align}

	Moreover, we have that $L(f, \Ad, 1) = L(1, \Ad^0 \varrho_0) \cdot L(1, \Ad^0 \varrho_0 \otimes \chi_F)$. Combining this with Equations~\eqref{eqn:<f_0,f_0>} and~\eqref{eqn:ratio_of_norms} gives:
	\begin{align}
		\langle f_0, f_0 \rangle & \sim_{E^\times} \frac{1}{f_{\Ad^0 \varrho}^{1/2}} \cdot \log|u_{f_0}|, \label{eqn:<f0,f0>=loguf0} \\
		\frac{\langle f, f \rangle}{\langle f_0, f_0 \rangle} & \sim_{E^\times} \frac{1}{f_{\Ad^0 \varrho}^{1/2}} \cdot \log|u_{f_0}^F|. \label{eqn:ratio=loguf0F}
	\end{align}
	Since \cite[Prop.\ 5.8]{Horawa} shows that $f_{\Ad^0 \varrho}$ is a square, this completes the proof.
\end{proof}

\subsection{Consequence for period invariants of Hilbert modular forms}\label{subsec:period_invariants}

Using Theorem~\ref{thm:motivic_action}, we are able to answer an old question raised by Michael Harris' work on period invariant of Hilbert modular forms~\cite{Harris_Periodinvariants}. For a Hilbert modular form $f$ of weight $(k_1, k_2; r)$ with $k_1 \equiv k_2 \equiv r \mod 2$, there are similar coherent cohomology classes $\omega_f^i \in H^1(X_\C, \mathcal E_i)_f$ for some rational line bundles $\mathcal E_i$. Under the assumption that $k_i \geq 2$, it turns out that $\dim H^1(X_\C, \mathcal E_i)_f = 1$, and Harris defines $\nu_f^i \in \C^\times$,  well-defined up to $E^\times$, as the constant such that:
$$\frac{\omega_f^i}{\nu_f^i} \in H^1(X, \mathcal E_i)_f.$$
These {\em period invariants} play a crucial role in the algebraicity of Rankin--Selberg and triple product $L$-functions of Hilbert modular forms. In special cases, they are related to CM periods and Shimura's periods~\cite{Harris_PeriodinvariantsII}; see \cite[Remark 4.11]{Horawa} for further discussion.

This raises the natural question: for $k_1 = k_2 = 1$, is any multiple of $\omega_f^i$ rational?

\begin{corollary}\label{cor:Harris_periods}
	Suppose $f$ is the base change of a weight one modular form $f_0$ of square-free level co-prime to $D$ whose Galois representation does not have projective image $A_5$. Then no multiple of $\omega_f^1$ or $\omega_f^2$ is $\overline \Q$-rational.
\end{corollary}
\begin{proof}
	Suppose that $\lambda \cdot \omega_f^1$ is $\overline \Q$-rational for some $\lambda \in \C$. Then $s^\ast(\lambda \cdot \omega_f^1) = \lambda \cdot \omega_f^2$ is rational. Therefore, $\lambda(\omega_f^1 + \omega_f^2),  \lambda(\omega_f^1 - \omega_f^2)\in H^1(X, \omega)_f \otimes_{\Q(f)} \overline \Q$, and comparing this with Theorem~\ref{thm:motivic_action} gives:
	$$\lambda \sim_{\overline \Q^\times } \log|u_{f_0}|^{-1}, \quad \lambda \sim_{\overline \Q^\times} \log|u_{f_0}^F|^{-1},$$
	and hence
	$$\log|u_{f_0}| \sim_{\overline \Q^\times} \log|u_{f_0}^F|.$$
	We claim this is impossible. Since the level of $f_0$ is co-prime to $D$, $\varrho_{f_0}$ is factors through a finite Galois extension $L/\Q$ linearly independent from $F/\Q$. Moreover, $\chi_F$ corresponds to the non-trivial map $\Gal(F/\Q) \to \{\pm 1\} \subseteq \Q^\times$. Consider a basis $u_1, \ldots, u_r$ of $\O_L^\times$ and extend it to a basis $u_1, \ldots, u_r, v_1, \ldots, v_s$ of $\O_{LF}^\times$. Then, writing the unit groups additively and $E$ for the field of definition of $\Ad^0 \varrho$, 
	\begin{align*}
		u_{f_0} & = \sum a_i u_i & u_i \in E^\times, \\
		u_{f_0}^F & = \sum b_i u_i + \sum c_j v_j & b_j, c_j \in E^\times.
	\end{align*}
	Writing $\sigma \in \Gal(F/\Q)$ for the non-trivial element, note that $\sigma u_f^F = - u_f^F$, while $\sigma u_i = u_i$, which shows that $c_j \neq 0$ for some $j$. On the other hand:
	\begin{align*}
		\log|u_{f_0}| & = \sum a_i \log|u_i| & u_i \in E^\times, \\
		\log|u_{f_0}^F| & = \sum b_i \log|u_i| + \sum c_j \log|v_j| & b_j, c_j \in E^\times.
	\end{align*}
	If $\mu \log|u_{f_0}| = \log |u_{f_0}^F|$ for some $\mu \in \overline\Q^\times$, then
	$$ \sum (\mu a_i - b_i) \log|u_i| - \sum c_j \log|v_j| = 0,$$
	which shows that $c_j = 0$ for all $j$ by Baker's linear independence of logarithms. This contradiction shows that $\log|u_{f_0}| \not\sim_{\overline \Q^\times} \log|u_{f_0}^F|$, and hence proves the corollary.
\end{proof}

\subsection{Explicit example}

In this section, we work out an explicit example of our result. More specifically, we will work out the identities~\eqref{eqn:<f0,f0>=loguf0} and~\eqref{eqn:ratio=loguf0F} to show how the first part of Theorem~\ref{thm:motivic_action} implies the second in a special case. This also gives a concrete example of Corollary~\ref{cor:Harris_periods}.

Let $F = \Q(\sqrt{D})$ where $D > 0$ is a fundamental discriminant. For~\eqref{eqn:<f0,f0>=loguf0}, we will use~\eqref{eqn:<f_0,f_0>} directly. For~\eqref{eqn:ratio=loguf0F}, we combine~\eqref{eqn:<f,f>} with~\eqref{eqn:<f_0,f_0>} to get the formula:
\begin{align}\label{eqn:ratio_of_norms}
	\frac{\langle f, f \rangle}{\langle f_0, f_0 \rangle} & \sim_{\Q^\times} \frac{|D|^{1/2}}{\pi^2} \frac{L(f, \Ad, 1)}{L(f_0, \Ad, 1)}.
\end{align}
Now, Artin formalism gives the factorization:
\begin{equation}
	L(f, \Ad, s) = L(f_0, \Ad, s) \cdot L(\Ad^0 \varrho_0 \otimes \chi_F, s),
\end{equation}
where $\chi_F$ is the quadratic character associated with $F/\Q$. Therefore,
\begin{equation}
	\frac{\langle f, f \rangle}{\langle f_0, f_0 \rangle}  \sim_{\Q^\times} \frac{|D|^{1/2}}{\pi^2}  L(\Ad^0 \varrho_0 \otimes \chi_F, 1).
\end{equation}
Therefore, computing an explicit version of~\eqref{eqn:ratio=loguf0F} amounts to finding a formula for $L(\Ad^0 \varrho_0 \otimes \chi_F, 1)$. 

We now specialize to the example of complex cubic fields as in \cite[\S 2.5]{Horawa}. Let $M/\Q$ be an imaginary quadratic extension with class number three, so that its class field $L$ is an $S_3$ extension of $\Q$. Let $K = L^{(123)}$ be the associated complex cubic field.

\begin{center}
	\begin{tikzcd}[scale = 2]
		& L \\
		M \ar[ru, no head, "3"]  & \ar[u, no head, swap, "2"]  K \\
		\Q \ar[u, no head, "2"] \ar[ru, no head, swap, "3"]
	\end{tikzcd}
\end{center}

We take $\varrho_0$ to be the unique irreducible representation of $S_3$, which gives rise to a weight one modular form $f_0$. Note that $\varrho_0 = \mathrm{Ind}_{M}^L \chi$ for the non-trivial character $\chi$ of $C_3$, and explicitly~$f = \theta_{\chi}$.

Then Stark observed that $u_{f_0} = \epsilon$, the generator of unit group of the cubic field $K$. 
Since $\varrho_0 = \Ind \chi$, $\Ad^0 \varrho_0 = \Ind(\chi) \oplus \chi_M$. Then one can check that:
\begin{equation}
	L(\Ad^0 \varrho_0, s) = L(\chi, s) \frac{\zeta_M(s)}{\zeta(s)} = \frac{\zeta_K(s) \zeta_M(s)}{\zeta(s)^2}.
\end{equation}
We recall the class number formula: for a number field $K$,
\begin{equation}
	\mathrm{Res}_{s=1} \zeta_K(s) = \frac{2^{r_1} \cdot (2 \pi)^{r_2} \cdot \mathrm{Reg}_K \cdot h_K}{w_K \cdot \sqrt{|D_K|}}.
\end{equation}
Then we obtain:
\begin{equation}\label{eqn:Stark_cubic_fd}
	L(\Ad^0 \varrho_0, 1)  =   \frac{6 \pi^2}{|D_M|} \log|\epsilon|.
\end{equation}

\begin{example}
	Let $M = \Q(\sqrt{-23})$. Then $\epsilon$ is a root of $x^3 - x - 1$. Combining~\eqref{eqn:<f_0,f_0>} with~\eqref{eqn:Stark_cubic_fd} gives:
	\begin{equation}
		\langle f_0, f_0 \rangle = \frac{1}{4\pi^2} \cdot \frac{23}{24} \cdot   \frac{6 \pi^2}{23} \log|\epsilon| = \frac{1}{16}  \log |\epsilon|.
	\end{equation}
	In Stark's classical normalization of the Petersson, this gives the identity~\cite[pp.\ 91]{Stark-II}
	$$\langle f_0, f_0 \rangle_{\mathrm{Stark}} = 3 \log|\epsilon|,$$
	because
	$$\langle f_0, f_0 \rangle_{\mathrm{Stark}} = 2 \cdot [{\rm PSL}_2(\mathbb Z) : \overline{\Gamma_0(23)}] \cdot \langle f_0, f_0 \rangle = 16 \cdot 3 \cdot \langle f_0, f_0 \rangle.$$
	(Here, we have included the rational constants which were left implicit in~\eqref{eqn:<f_0,f_0>} to recover Stark's result on the dot.) 
	
	Alternatively, in this example, we could have used the formula
	$$L(\Ad^0\varrho_0,s) = \frac{\zeta_L(s)}{\zeta_K(s)}$$
	and observe that $\Reg_L = 3 \Reg_K^2 = 3 \log|\epsilon|^2$. Indeed, the unit group of $L$ is generated by $\epsilon \in \mathcal O_K^\times$ and a unit $\delta \in \mathcal O_L^\times$ such that $N_{L/K} \delta = \epsilon$. For embeddings $\sigma_1, \sigma_2, \sigma_3 \colon L \hookrightarrow \C$, the regulator is:
	\begin{align*}
		\det \begin{pmatrix}
			2\log|\sigma_1(\epsilon)| & 2\log |\sigma_2(\epsilon)| \\
			2\log|\sigma_1(\delta)| & 2\log |\sigma_2(\delta)|
		\end{pmatrix} & = 2 \det \begin{pmatrix}
			2 \log |\epsilon | & - \log|\epsilon| \\
			\log|\epsilon| & \log|\epsilon|
		\end{pmatrix} \\
		& = 3 \log|\epsilon|^2.
	\end{align*}
\end{example}

Next, we would like to understand the $L$-function $L(\Ad^0 \varrho_0 \otimes \chi_F, s)$ in terms of Dedekind $\zeta$-functions. We have a bigger diagram of fields:

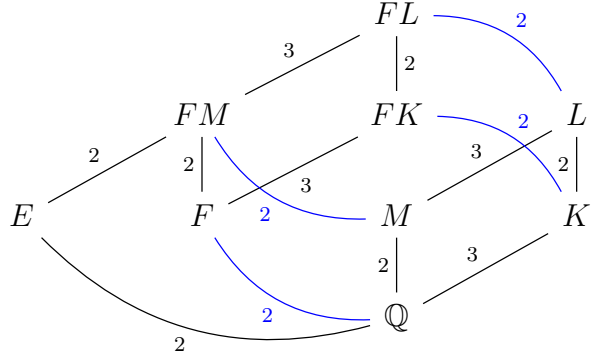
\begin{center}
	\begin{tikzcd}[column sep = large]
		& & FL \\
		& FM \ar[ru, no head, "3"]  & \ar[u, no head, swap, "2"]  FK & \ar[lu, no head, bend right, swap, blue, "2"] L \\
		E \ar[ru, no head, "2"] & 	F \ar[u, no head, "2"] \ar[ru, no head, swap, "3"] & M \ar[lu, no head, bend left, blue, "2"] \ar[ru, no head, "3"] & K \ar[lu, no head, bend right, swap, blue, "2"] \ar[u, no head, "2"] \\
		& &  \ar[llu, bend left, no head, "2"]\Q \ar[lu, no head, bend left, blue, "2"] \ar[u, no head, "2"] \ar[ru, no head, "3"] \\
	\end{tikzcd}
\end{center}
where $E$ is the other intermediate field between $FM$ and $\Q$. 

Recalling that $\Ad^0 \varrho_0 = \varrho_0 \oplus \chi_M$, we have $\Ad^0 \varrho_0 \otimes \chi_F = \varrho_0 \otimes \chi_F \oplus \chi_E$. One then check that:
$$L(\Ad^0 \varrho_0 \otimes \chi_F, s) = L(\varrho_0 \otimes \chi_F, s) \cdot \frac{\zeta_E(s)}{\zeta(s)}.$$
To get a formula for $L(\varrho_0 \otimes \chi_F, s)$, we observe that
$$\frac{\zeta_{FL}(s)}{\zeta_K(s) \cdot \zeta_{FK}(s)} = \frac{\zeta_M(s) \cdot L(\varrho_0 \otimes \chi_F, s) \cdot \zeta_E(s)}{\zeta(s)^3}$$
and rearrange this to get
\begin{equation}
	L(\Ad^0 \varrho_0 \otimes \chi_F, s) =  \frac{\zeta_{FL}(s)}{\zeta_K(s) \cdot \zeta_{FK}(s)} \cdot \frac{\zeta(s)^3}{\zeta_M(s) \zeta_E(s)} \cdot \frac{\zeta_E(s)}{\zeta(s)} = \frac{\zeta_{FL}(s) \zeta(s)^2}{\zeta_K(s) \zeta_{FK}(s) \zeta_M(s) }.
\end{equation}

In fact, there is a simpler formula. Since:
\begin{align*}
	\zeta_K(s) & = L(\varrho_0, s) \cdot \zeta(s), \\
	\zeta_{KF}(s) & = L(\varrho_0 \otimes \chi_F, s) \cdot L(\varrho_0, s) \cdot \zeta_F(s),
\end{align*}
we see that
\begin{equation}\label{eqn:Ad0_twist}
	L(\Ad^0 \varrho_0 \otimes \chi_F, s) = L(\varrho_0 \otimes \chi_F, s) \cdot \frac{\zeta_E(s)}{\zeta(s)} = \frac{\zeta_{FK}(s) \zeta(s)}{\zeta_K(s) \zeta_F(s)} \cdot \frac{\zeta_E(s)}{\zeta(s)} = \frac{\zeta_{FK}(s)  \zeta_E(s)}{\zeta_K(s) \zeta_F(s)}.
\end{equation}

\begin{example}
	Let $M = \Q(\sqrt{-23})$ and $F = \Q(\sqrt{5})$. We record the relevant quantities:
	$$\begin{array}{c|c|c|c|c|c}
		& r_1 & r_2 & h & w & |D| \\ \hline
		FK & 2 & 2 & 1 & 2 & 5^3 \cdot 23^2   \\ \hline
		K & 1 & 1 & 1 & 2 & 23  \\ \hline
		F & 2 & 0 & 1 & 2 & 5 \\ \hline
		E &  0 & 1 & 2 & 2 & 5 \cdot 23
	\end{array}$$
	Then by Equation~\eqref{eqn:Ad0_twist}, we get
	\begin{equation}
		L(\Ad^0 \varrho_0 \otimes \chi_F, 1) = \frac{\frac{2^{2} \cdot (2 \pi)^{2} \cdot 1}{2 \cdot 5 \cdot \sqrt{5} \cdot 23}  \cdot \frac{2^{0} \cdot (2 \pi)^{1} \cdot 2}{2 \cdot \sqrt{5 \cdot 23}} }{\frac{2^{1} \cdot (2 \pi)^{1} \cdot 1}{2 \cdot \sqrt{23}} \cdot \frac{2^{2} \cdot (2 \pi)^{0} \cdot 1}{2 \cdot \sqrt{5}}} \cdot \frac{\Reg_{FK} }{\Reg_K \Reg_F} = \frac{(2\pi)^2}{5 \cdot \sqrt{5} \cdot 23} \cdot \frac{\Reg_{FK}}{\Reg_K \Reg_F}
	\end{equation}
	
	The unit group of $FK$ is generated by three units: $\epsilon_K$, $\epsilon_F$ and a third unit $\delta$ such that $N_{KF/K} \delta = -1$. For three embeddings $\sigma_1, \sigma_2, \sigma_3 \colon FK \hookrightarrow \C$, the regulator is
	\begin{align*}
		\Reg_{FK} & = \det \begin{pmatrix}
			\log|\sigma_1(\epsilon_K)| & \log|\sigma_2(\epsilon_K)| & 2\log|\sigma_3(\epsilon_K)|  \\
			\log|\sigma_1(\epsilon_F)| & \log|\sigma_2(\epsilon_F)| & 2\log|\sigma_3(\epsilon_F)| \\
			\log|\sigma_1(\delta)| & \log|\sigma_2(\delta)| & 2\log|\sigma_3(\delta)|
		\end{pmatrix} \\
		& = \det \begin{pmatrix}
			\log|\epsilon_K| & \log|\epsilon_K| & -\log|\epsilon_K|  \\
			\log|\epsilon_F| & - \log|\epsilon_F| & 2\log|\epsilon_F| \\
			\log|\sigma_1(\delta)| & -\log|\sigma_1(\delta)| & 2\log|\sigma_3(\delta)| \\
		\end{pmatrix} \\
		& = \log|\epsilon_K| \log|\epsilon_F| \det \begin{pmatrix}
			1 & 1 & -1 \\
			1 & -1 & 2 \\
			\log|\sigma_1(\delta)| & -\log|\sigma_1(\delta)| & 2\log|\sigma_3(\delta)| \\
		\end{pmatrix} \\
		& = 4 \Reg_K \cdot  \Reg_F \cdot (\log|\sigma_1(\delta)/\sigma_3(\delta)|)
	\end{align*}
	
	Finally, we get from~\eqref{eqn:ratio_of_norms} that:
	\begin{equation}\label{eqn:example_of_ratio}
		\frac{\langle f, f \rangle}{\langle f_0, f_0 \rangle} \sim_{\Q^\times} \frac{\sqrt{D}}{\pi^2} \frac{4\pi^2}{5 \sqrt{5} \cdot 23} \frac{\Reg_{FK}}{\Reg_K \Reg_F} \sim_{\Q^\times}  \frac{\Reg_{FK}}{\Reg_K \Reg_F} \sim_{\Q^\times} \log\left|\frac{\sigma_1(\delta)}{\sigma_3(\delta)}\right|.
	\end{equation}
	
\end{example}

 \printbibliography
\end{document}